\documentclass[italian,a4paper,11pt]{article}
\usepackage[latin1]{inputenc}
\usepackage[english]{babel}
\usepackage{babel,amsmath,amssymb,amsbsy,amsfonts,latexsym,amsthm,mathrsfs}

\input pictex.sty
\input epsf.tex

\usepackage{amsmath,amsfonts,latexsym}
\newtheorem{theorem}{Theorem}[section]
\newtheorem{proposition}[theorem]{Proposition}

\newtheorem{lemma}[theorem]{Lemma}

\newtheorem{remark}[theorem]{Remark}
\newtheorem{example}[theorem]{Example}
\newtheorem{assumption}[theorem]{Assumption}
\newtheorem{hypothesis}[theorem]{Hypothesis}
\def\proof{{\sl Proof. }}
\def\cl#1{{\mathscr #1}}

\def\P{{\mathbb P}}
\def\R{{\mathbb R}}

\def\N{{\mathbb  N}}

\def\E{{\mathbb E}}

\def\Var{{\rm Var }}
\def\Cov{{\rm Cov }}

\def\<{\langle}
\def\>{\rangle}

\newcommand{\cvd}{
                   $\quad\Box $
                   \medskip
                   }

\def\appendix{\par
  \setcounter{chapter}{0}
  \def\@chapter{Appendix}
  \def\thechapter{\Alph{chapter}}}
\font\tengoth=eufm10 scaled 1200 \font\sevengoth=eufm7 scaled 1200
\font\fivegoth=eufm5 scaled 1200
\newfam\gothfam
\textfont\gothfam=\tengoth \scriptfont\gothfam=\sevengoth
\scriptscriptfont\gothfam=\fivegoth

\def\coeffbin#1#2{{\Big(\begin{array}{c}#1\\ #2\end{array}\Big)}}

\def\ep{\varepsilon}

\begin{document}

\title{
\bf
Large deviation estimates of the crossing probability
for pinned Gaussian processes}
\author{
{\sc Lucia Caramellino}\smallskip\\
{\sc Barbara Pacchiarotti}\smallskip\\
 Dept. of Mathematics, University of  Rome-Tor Vergata \smallskip\\
{\small mailto:
\texttt{\{caramell,pacchiar\}@mat.uniroma2.it}}}
%
\maketitle

{\small
\bf Abstract. \rm The paper deals with the asymptotic behavior of
the bridge of a Gaussian process conditioned to
stay in $n$ fixed points at $n$ fixed past instants.
In particular, functional large deviation results are stated
for small time. Several examples
are considered: integrated or not fractional Brownian motion,
$m$-fold integrated Brownian motion.
As an application, the asymptotic behavior of the exit
probability is studied and used
for the practical purpose of the numerical computation, via
Monte Carlo methods, of
the hitting probability up to a given time.
}

\vskip 1cm


\bigskip


\vfill
\bf Corresponding Author\rm: Lucia Caramellino,
Dipartimento di Matematica,
Universit\`a di Roma \sl Tor Vergata\rm, Via della Ricerca
Scientifica, I-00133 Roma, Italy; fax: + 39 06 72594699.

\section{Introduction}
Simple formulas for the crossing probability in small time
for pinned processes have
been recently investigated in the literature, because
of their use in improving the performance of the numerical simulation
of processes to be killed when a prescribed boundary is reached.
The idea underlying the application is simple. In fact, consider the generical
step of the simulation procedure: one has generated the process
of interest, say $X$, at some $n\geq  1$ fixed instants $0<T_1<\cdots<T_n$,
observing the positions $X_{T_1}=x_1,\ldots, X_{T_n}=x_n$.
In order to get the exit time,
one simulates again the process, at time $T_n+\ep$, and
if the observed position $X_{T_n+\ep}=y$ reaches the boundary,
then the crossing is achieved.
This gives rise to an over estimate  of the exit time, 
which can dramatically bring to a significant error, as observed by
many authors. One way to overcome this difficulty is to compute the
crossing probability of the pinned process, that is for the
conditional process
$(X_{T_n+\ep t})_{0\leq t\leq 1}$ given all the past
observations $X_{T_1}=x_1,\ldots, X_{T_n}=x_n$
and the present one $X_{T_n+\ep}=y$, and to use it in order to decide if the boundary has been
breached or not. Let us stress that in the general case, no closed formulas
are available, so that such a procedure is carried out with an approximation
(by large deviations, as $\ep\to 0$) of the exit probability.

In the case of diffusion processes, the Markov property
allows one to work with the bridge-process between the observations at times
$T_n$ and $T_n+\ep$. This case has bee widely studied in the
literature, see e.g. Baldi
and Caramellino \cite{bib:bc} and the references quoted therein.
This approach obviously fails if a non Markovian process
is taken into account, so that one has to consider all the past
observations and to handle the bridge of the conditional process.

The present paper deals with the large deviation asymptotic behavior of the exit probability
of such a pinned process whenever the original one is a (continuous)
Gaussian process (and in particular, not necessarily a Markovian one).
Our wide class of examples can be split in two main sets.

First, we consider the fractional
Brownian motion, which  is widely used in risk theory modelling
(see e.g. Baldi and Pacchiarotti \cite{bib:bp}). As a consequence,
 we can handle the semimartingale process
resulting from a linear combination between a fractional Brownian motion
with Hurst index greater than 3/4 and
a standard Brownian motion, independent each other,
a promising tool to set up a
non Markovian model in Mathematical Finance  (see Cheridito \cite{bib:cheridito1}).

Secondly, we can deal with an integrated Gaussian process,
that is a process defined as
the integral w.r.t. the Lebesgue measure
of a Gaussian one. As an example, we obtain the integrated fractional Brownian motion,
which is linked to fractal
properties of solutions of the inviscid
Burgers equation. Notice that the law of its maximum, analyzed  e.g.
in Molchan and Khokhlov \cite{bib:mk}, is strictly connected
to level crossing probabilities.
Furthermore, we can consider  the $m$-fold iterated Brownian motion (see e.g.
Chen and Li \cite{bib:cl}) and in particular  the
integrated Brownian
motion, having interesting applications in nonparametric
estimating in Statistics (see e.g. Groeneboom, Jongbloed and Wellner
\cite{bib:gjw} and references quoted therein) and
used in metrology as a model for the atomic clock error
(whose precision and re-synchronization are strictly
related to the level crossing, see e.g. Galleani, Sacerdote, Tavella and Zucca
\cite{bib:torino}).

The paper is organized as follows.
After a brief recall of some well known results related to large deviations
for Gaussian processes (Section \ref{tutorial}), we first
get a functional large deviation result, for small
time, for Gaussian processes conditioned
to stay in $n$ fixed positions $x_1,\ldots, x_n$ at $n$ fixed
instants $T_1<\cdots<T_n$ (Section \ref{model}). In a second
moment (Section \ref{bridge}), we state a functional large
deviation principle for the bridge of such conditional processes.
Let us stress that, surprisingly, we obtain a degenerate kind of
large deviations for Gaussian processes
having a quite smooth covariance function (e.g. for
integrated Gaussian processes), so that we give some refined results
allowing to handle also these cases. In particular, we obtain examples
of an interesting and non trivial asymptotic behavior, in
which the (non degenerate) large deviation speed is different according to
the conditional process or its bridge: the speed associated to the
bridge can be much faster than the one turning out for
the conditional process.
Finally (Section \ref{exit}),
we are able to give the asymptotic behavior, in terms of large
deviations, of the probability of crossing one or two
possibly time dependent levels, and we propose some numerical results
concerning the fractional Brownian motion.

\section{Large deviations for Gaussian processes}\label{tutorial}

We briefly recall here some main facts related to the
large deviation theory for Gaussian processes we are going to use.
There are many references in the literature on this topics, where
all details and proofs may be found; let us here
recall some classical references: Azencott \cite{bib:azencott},
Deuschel and Strook \cite{bib:ds}, Dembo and Zeitouni
\cite{bib:dz}. Without loss
of generality, we can consider centered Gaussian processes.

Throughout the paper, $C([0,1])$ will denote
the set of the continuous paths on $[0,1]$,
endowed with the topology induced by the sup-norm.
Moreover, $\cl M[0,1]$ will be its dual, i.e.
the set the signed Borel measures
on $[0,1]$, and for any $\lambda\in \cl M[0,1]$,
$\<\lambda,\cdot\>$ will stand for the associated linear functional:
$\<\lambda,h\>=\int_0^1h_t\lambda(dt)$, $h\in C([0,1])$.

A continuous process $U=(U_t)_{t\in [0,1]}$,
defined on some probability space $(\Omega,\cl F, \P)$,
is a centered Gaussian
process if for any $\lambda\in \cl M[0,1]$ then
$\<\lambda, U\>=\int_0^1U_t\lambda(dt)$ is a centered Gaussian
r.v. taking value on $\R$. The associated
continuous covariance function
$k(t,s)=\Cov(U_t, U_s)$,  $t,s\in[0,1]$,
plays  a crucial role. For example, one has
$$
\Var(\<\lambda, U\>)=\int_0^1\int_0^1k(t,s)\lambda(dt)\lambda(ds),
\quad\mbox{for any }\lambda\in\cl M[0,1].
$$
In addition to $k$, another important instrument for handling
Gaussian processes is the associated reproducing kernel Hilbert
space. It is a Hilbert space in $C([0,1])$ which is usually defined through the following dense
subset:
$$
\cl L=\Big\{h\in C([0,1])\,:\, h_t=\int_0^1k(t,s)\lambda(ds), \mbox{
with } \lambda\in\cl M[0,1]\Big\}.
$$
Let us be a little bit
more precise about $\cl H$. First, let $\mu$ denote the
measure induced by the Gaussian process $U$:
$\mu(A)=\P(U\in A)$ for any Borel set $A$ in $C([0,1])$.
Let $\Gamma\subset L^2(\mu)$ be defined as the following set of
(real) Gaussian r.v.'s:
$$
\Gamma=\{Y\,:\,Y(\cdot)=\<\lambda,\cdot\>, \mbox{ with }
\lambda\in\cl M[0,1]\}.
$$
It immediately follows that for $Y_1,Y_2\in \Gamma$,
with $Y_i(\cdot)=\<\lambda_i,\cdot\>$ as $i=1,2$, then
$$
\Cov(Y_1,Y_2)=(Y_1,Y_2)_{L^2(\mu)}=\int_0^1\int_0^1
k(t,s)\lambda_1(dt)\lambda_2(ds),
$$
where, from now on, the symbol $(\cdot,\cdot)$ denotes
an inner product. We define now
$$
H={\overline{\phantom{|}\Gamma}}^{\ \|\cdot\|_{L^2(\mu)}}.
$$
Obviously, $H$ is a  closed subspace of $L^2(\mu)$
and is indeed a set
of Gaussian r.v.'s taking values on $\R$. Moreover, it becomes a Hilbert space
if endowed with the inner product
$$
(Y_1,Y_2)_H=(Y_1,Y_2)_{L^2(\mu)},\quad Y_1,Y_2\in H.
$$
We now set the following map:
$$
\begin{array}{lll}
{\cal S} : & H \rightarrow & C([0,1]) \\
    & Y \mapsto &(SY)_t=\displaystyle\int x_t Y(x)\mu(dx)\equiv \E(U_t Y).
\end{array}
$$
It can be shown that ${\cal S}$ is a linear, one-to-one and continuous
map, so that ${\cal S}^{-1}\,:\, {\cal S} H\to H$ is a well posed
continuous and linear map.
The reproducing kernel Hilbert space $\cl H$ is defined as the
image  of $H$ through ${\cal S}$:
$$
\cl H={\cal S}H\equiv\{h\in C([0,1])\,:\, h_t=(SY)_t, Y\in H\}.
$$
Finally, setting
$$
(h_1,h_2)_{\cl H}=({\cal S}^{-1}h_1, {\cal S}^{-1}h_2)_H
\equiv({\cal S}^{-1}h_1,
{\cal S}^{-1}h_2)_{L^2(\mu)},\quad h_1,h_2\in \cl H,
$$
then $(\cdot, \cdot)_{\cl H}$ is an inner product on $\cl H$,
which in turn makes $\cl H$ a Hilbert space. This is the rigorous
definition of the reproducing kernel Hilbert space associated to a
(centered) Gaussian process.  Finally, it immediately follows that
$$
\cl H={\overline{\cl L}}^{\ \|\cdot\|_{\cl H}},
\quad\mbox{with } \cl L=\{x\in C([0,1])\,:\, x_t=\int_0^1
k(t,s)\lambda(ds),\mbox{ with }\lambda\in \cl M[0,1]\}.
$$

In the sequel, we will speak about ``the reproducing kernel
Hilbert space associated to the covariance function $k(t,s)$''.
In fact, given a continuous, symmetric and positive
definite  function $k(t,s)$ defined on
$[0,1]\times[0,1]$, one can build a centered and
continuous Gaussian process
$U=(U_t)_{t\in [0,1]}$ having $k$ as its covariance function.
Now, the associated reproducing kernel Hilbert space
is naturally defined.

The main property we are going to use
is related to the Cram\`er transform:
\begin{theorem}\label{cramer-transform}{\rm\bf[Cram\`er transform]}
Let $I$ denote the Cram\`er transform, that is
$$
I(x)=\sup_{\lambda\in\cl M[0,1]}\Big(\<\lambda , x\>-\log
\E(e^{\<\lambda, U\>})\Big)
=\sup_{\lambda\in\cl M[0,1]}\Big(\<\lambda , x\>-\frac 12
\int_0^1\int_0^1 k(t,s)\lambda(dt)\lambda(ds)\Big).
$$
Then,
$$
I(x)=\left\{
\begin{array}{ll}
\displaystyle\frac 12 \,\|x\|_{\cl H}^2 & \mbox{ if } x\in\cl H\smallskip\\
+\infty & \mbox{ otherwise.}
\end{array}
\right.
$$

\end{theorem}

Suppose now to have a family of continuous Gaussian
processes $\{U^\ep\}_\ep$: is it possible to determine a
large deviation principle? Because of
the special form of the Laplace transform for Gaussian measures,
a large deviation principle can be stated if a nice asymptotic
behavior holds for the Laplace transforms, as summarized
in the following
\begin{theorem}\label{GD-G}
Let $\{U^\ep\}_\ep$ be a family of continuous Gaussian processes.
Let $\gamma_\ep$ be an infinitesimal function, i.e. $\lim_{\ep\to 0}
\gamma_\ep=0$, and suppose that, for any $\lambda\in \cl M[0,1]$,
$$
0=\lim_{\ep\to 0}\E(\<\lambda,U^\ep\>)
\quad{and}\quad
\Lambda(\lambda)=\lim_{\ep\to 0}
\displaystyle\frac {\Var(\<\lambda, U^\ep\>)}{\gamma^2_\ep}
\equiv\int_0^1\int_0^1 \bar k(t,s)\lambda(dt)\lambda(ds),
$$
for some continuous, symmetric and positive definite function
$\bar k$.
Then, $\{U^\ep\}_\ep$ satisfies a large deviation principle
on $C([0,1])$,
with inverse speed $\gamma^2_\ep$ and (good) rate function
\begin{equation}\label{RF-G}
I(h)=\left\{
\begin{array}{cl}
\displaystyle\frac 12 \|h\|^2_{\bar{\cl H}}
& \mbox{if }h\in\bar{\cl H}\smallskip\\
+\infty & \mbox{otherwise},
\end{array}
\right.
\end{equation}
where $\bar{\cl H}$ and $\|\cdot\|_{\bar{\cl H}}$ denote, respectively,
the reproducing kernel Hilbert space and the related norm
associated to the covariance function $\bar k$.
\end{theorem}
Let us recall, once for all, that the sentence
``$\{U^\ep\}_\ep$ satisfies a large deviation principle on $C([0,1])$
with inverse speed $\gamma^2_\ep$ and (good) rate function $I$''
means: $\lim_{\ep\to 0}\gamma_\ep=0$; the set
$\{I\leq a\}$ is compact in $C([0,1])$, for any fixed $a$;
 the following inequalities hold:
\begin{itemize}
\item[-] for any open set $G$ in $ C([0,1])$,
$\liminf_{\ep\to 0}\gamma^2_\ep\log \P(U^\ep\in
G)\geq -\inf_{h\in G} I(h)$;
\item[-] for any closed set $F$ in $C([0,1])$,
$\limsup_{\ep\to 0}\gamma^2_\ep
\log \P(U^\ep\in F)\leq -\inf_{h\in F} I(h)$.
\end{itemize}

For the sake of convenience,  Theorem \ref{GD-G} is written for
a non-centered family of Gaussian processes, even if it requires that
the expected path weakly converges to zero.
The idea of the proof of Theorem \ref{GD-G} if the following.
It is well known (e.g. by applying
the Gartner-Ellis Theorem, see e.g. Dembo and Zeitouni \cite{bib:dz})
that a large deviation principle holds if the hypotheses of
Theorem \ref{GD-G} are satisfied, and the rate function
is given by the Legendre transform of
$$
\bar \Lambda(\lambda)
=\frac 12\, \int_0^1\int_0^1 \bar k(t,s)\lambda(dt)\lambda(ds),\quad
\lambda\in \cl M[0,1].
$$
In view of Theorem \ref{cramer-transform}, one immediately
obtains formula (\ref{RF-G}).

\section{Large deviations for the conditional process}\label{model}
Let $X=(X_t)_{t\ge 0}$  be a  Gaussian, centered  process with
continuous covariance function
\begin{equation}\label{k}
k(t,s)=\Cov(X_t,X_s).
\end{equation}
For a fixed $n\in \N$ and $j=1,\ldots, n$, let $X^j=(X^j_t)_{t\ge 0}$
stand for the process giving the conditional behavior of $X$ given that
it  assumes the values $x_1,\ldots,x_j$ at the
$j$ times $0<T_j<\ldots<T_j$  respectively. Since the original process
$X$ is Gaussian, the process
$X^j=(X^j_t)_{t\ge 0}$  is equal in law to
\begin{equation}\label{cond-process}
X^j_t=X^{j-1}_t - \alpha_{j}(t)(X^{j-1}_{T_j}-x_j),
\end{equation}
where
\begin{equation}\label{alphaj}
\alpha_j(t)=\frac {k_{j-1}(t,T_j)}{k_{j-1}(T_j,T_j)}
\end{equation}
and also $k_j$, giving the covariance function
associated to $X^j$, is recursively defined as
\begin{equation} \label{kj}
\begin{array}{rl}
k_j(t,s)=\Cov(X^j_t,X^j_s)=&k_{j-1}(t,s)-\alpha_j(t) k_{j-1}(s,T_j)\\
=&k_{j-1}(t,s)-\alpha_j(s) k_{j-1}(t,T_j).
\end{array}
\end{equation}
Obviously, the case $j=0$ is related to the original process and its
covariance function, that is $X^0\equiv X$ and $k_0\equiv k$.

Our first aim is to study the behavior of the covariance function
of the original process $X$ in order to get a functional large
deviation principle for the $n$-fold conditional process $X^n$
for small time, that is for $\{X^n_{T_n+\ep\cdot}\}_\ep$ as
$\ep\to 0$.

Let us consider an infinitesimal function $\gamma_\ep$ ($\gamma_\ep
\to 0$ as $\ep\to 0$), whose square will play the role of the inverse
speed of the large deviation principles we are going to study.

\begin{assumption}\label{ass-GD}
There exists the asymptotic covariance function $\bar
k(t,s)$, defined as
\begin{equation}\label{kbar}
\begin{array}{rl}
\bar k(t,s) =&\displaystyle\lim_{\ep\to 0} \frac {\Cov(X_{T_n+\ep
t}-X_{T_n},
X_{T_n+\ep s}-X_{T_n})}{\gamma^2_\ep}\\
=&\displaystyle\lim_{\ep\to 0} \frac {k(T_n+\ep t, T_n+\ep s)
-k(T_n+\ep t, T_n)-k(T_n, T_n+\ep s) +k(T_n, T_n)}{\gamma^2_\ep}
\end{array}
\end{equation}
uniformly as $(t,s)\in[0,1]\times[0,1]$.

\end{assumption}

As an immediate application of Theorem \ref{GD-G}
(take $U^\ep_t=X_{T_n+\ep t}-X_{T_n}$),
Assumption \ref{ass-GD} implies that
the family $\{(X_{T_n+\ep
t}-X_{T_n})_{t\in[0,1]}\}_\ep$ satisfies a large deviation
principle on $C([0,1])$, with inverse speed $\gamma^2_\ep$ and
good rate function given by
\begin{equation}\label{J}
J(h)=\left\{
\begin{array}{cl}
\displaystyle\frac 12 \|h\|^2_{\bar{\cl H}}
& \mbox{if }h\in\bar{\cl H}\smallskip\\
+\infty & \mbox{otherwise}
\end{array}
\right.
\end{equation}
where $\bar{\cl H}$ is the reproducing kernel Hilbert space
associated to the  covariance function $\bar k(t,s)$ as
and the symbol $\|\cdot \|_{\bar{\cl H}}$
denotes the usual norm defined on $\bar{\cl H}$.

Now, in order to achieve a large deviation principle for the
$n$-fold conditional process $X^n$, we have to investigate the
behavior of the functions $k_j$, defined through (\ref{kj}), in a
small time interval of length $\ep$. Let us consider the following

\begin{assumption}\label{ass-A}
For any fixed $T>0$, the following
limit exists:
\begin{equation}\label{rhobar}
\bar \rho(t,T)=\lim_{\ep\to 0}\frac{k(T_n+\ep
t,T)-k(T_n,T)}{\gamma_\ep}, \quad\mbox{ uniformly as }t\in[0,1].
\end{equation}
\end{assumption}

Let us discuss two simple
but useful  consequences of the assumptions introduced above.
\begin{lemma}\label{bar}
\begin{itemize}

\item[(i)]
Under Assumption \ref{ass-A}, as $j=1,\ldots,n$ one has
$$
\displaystyle\lim_{\ep\to
0}\frac{\alpha_j(T_n+\ep t) -\alpha_j(T_n)}{\gamma_\ep}=\bar \alpha_j(t), \quad
\mbox{uniformly as $t\in[0,1]$},
$$
where
\begin{equation}\label{alphajbar}
\bar\alpha_j(t)=\frac{\bar\rho_{j-1}(t,T_j)}{k_{j-1}(T_j,T_j)}
\end{equation}
$k_{j-1}$ being defined in (\ref{kj}), $\bar \rho_0\equiv \bar\rho$ and
\begin{equation}\label{rhojbar}
\begin{array}{rl}
\bar\rho_j(t,T)= \displaystyle\lim_{\ep\to 0}\frac{k_j(T_n+\ep
t,T)-k_j(T_n,T)}{\gamma_\ep}
=&\bar\rho_{j-1}(t,T)-\bar\alpha_j(t)k_{j-1}(T,T_j)\\
=&\bar\rho_{j-1}(t,T)-\alpha_j(T)\bar\rho_{j-1}(t,T_j),
\end{array}
\end{equation}
the above limit being uniformly as $t\in[0,1]$.

\item[(ii)]
Under Assumptions \ref{ass-GD} and \ref{ass-A},
as $j=1,\ldots,n$ one has
\begin{equation}\label{mjbar}
\lim_{\ep\to 0} \E(X^j_{T_n+\ep
t}-X^j_{T_n})=0, \quad \mbox{uniformly as } t\in[0,1]
\end{equation}
and
$$
\displaystyle\lim_{\ep\to 0} \frac
{\Cov(X^j_{T_n+\ep t}-X^j_{T_n}, X^j_{T_n+\ep
s}-X^j_{T_n})}{\gamma^2_\ep}=\bar k_j(t,s),\quad \mbox{uniformly as } t,s\in[0,1],
$$
with
\begin{equation}\label{kjbar}
\bar k_j(t,s)=\bar k(t,s)-\displaystyle\sum_{\ell=1}^j
k_{\ell-1}(T_\ell,T_\ell)\bar \alpha_\ell(t)\bar \alpha_\ell(s),
\end{equation}
with $\bar\alpha_\ell$ defined in (\ref{alphajbar}).
\end{itemize}
\end{lemma}

\proof
$(i)$
From Assumption \ref{ass-A} and  (\ref{alphaj}), one immediately
has
$$
\bar\alpha_1(t)=\frac{\bar\rho(t,T_1)}{k(T_1,T_1)}.
$$
Therefore, by using (\ref{kj}),
there exists, uniformly as $t\in[0,1]$,
$$
\begin{array}{rl}
\bar\rho_1(t,T)= \displaystyle\lim_{\ep\to 0}\frac{k(T_n+\ep
t,T)-k(T_n,T)}{\gamma_\ep}
=&\bar\rho_0(t,T)-\bar\alpha_1(t)k_0(T,T_1)\\
=&\bar\rho_0(t,T)-\alpha_1(T)\bar\rho_0(t,T_1)
\end{array}
$$
where, as usual, we have set $\bar\rho_0\equiv\bar\rho$ and
$k_0\equiv k$. This ensures the existence of $\bar\alpha_2$.
The statement now follows by iteration.

\smallskip

 The proof of $(ii)$ is a straightforward application of Assumption
\ref{ass-GD} and part $(i)$. \cvd

Notice that in particular, since $X^n_{T_n}=x_n$, one has,
again uniformly for $t,s\in[0,1]$,
\begin{equation}\label{mnbar}
x_n=\lim_{\ep\to 0} \E(X^n_{T_n+\ep t})
\end{equation}
\begin{equation}\label{knbar}
\bar k_n(t,s)=\lim_{\ep\to 0} \frac {\Cov(X^n_{T_n+\ep t},
X^n_{T_n+\ep s})}{\gamma^2_\ep}=\bar
k(t,s)-\displaystyle\sum_{\ell=1}^n k_{\ell-1}(T_\ell,T_\ell)\bar
\alpha_\ell(t)\bar \alpha_\ell(s).
\end{equation}

\medskip

We are now ready to prove the main large deviation result of this
section:

\begin{theorem}\label{th-AB}
Under Assumption \ref{ass-GD} and \ref{ass-A},
 the family $\{(X^n_{T_n+\ep
t})_{t\in[0,1]}\}_\ep$ satisfies a large deviation
principle on $C([0,1])$, with inverse speed $\gamma^2_\ep$ and
good rate function
\begin{equation}\label{Jn}
J_n(h)=\left\{
\begin{array}{cl}
\displaystyle \frac 12\, \|h-x_n\|^2_{\bar{\cl H}_n} & \mbox{ if }
h_0=x_n
\quad \mbox{and}\quad h-x_n\in \bar{\cl H}_n\\
+\infty & \rm{ otherwise.}
\end{array}
\right.
\end{equation}
$\bar{\cl H}_n$ being the reproducing kernel Hilbert space
associated to the covariance function
\begin{equation}\label{kbar-n}
\bar k_n(t,s)=\bar k(t,s)-\sum_{j=1}^n k_{j-1}(T_j,T_j)\bar
\alpha_j(t)\bar\alpha_j(s)
\end{equation}
where $\bar k(\cdot,\cdot)$, $k_j(\cdot,\cdot)$ and
$\bar\alpha_j(\cdot)$ are defined through (\ref{kbar}), (\ref{kj})
and (\ref{alphajbar}) respectively.
\end{theorem}

\proof
We start by showing that $\{(X^1_{T_n+\ep
t}-X^1_{T_n})_{t\in[0,1]}\}_\ep$, $X^1$ being defined in
(\ref{cond-process}) with $j=1$, satisfies a large deviation
principle.
By (\ref{mjbar}), it follows that
$$
\lim_{\ep\to 0} \E(\langle \lambda, X^1_{T_n+\ep \cdot}-X^1_{T_n}\rangle)
=\lim_{\ep\to 0}\int_0^1 \E( X^1_{T_n+\ep t}-X^1_{T_n})\lambda(dt)
=0
$$
and, recalling that $\Var(X_{T_1})=k(T_1,T_1)\equiv
k_0(T_1,T_1)$, by (\ref{kjbar})
$$
\displaylines{
\lim_{\ep\to 0} \frac{\Var(\langle \lambda, X^1_{T_n+\ep
\cdot}-X^1_{T_n}\rangle)}{\gamma^2_\ep}
=\lim_{\ep\to 0}
\int_0^1\!\!\lambda(dt)\!\int_0^1\lambda(ds)
\frac{\Cov(X^j_{T_n+\ep t}-X^j_{T_n}, X^j_{T_n+\ep
s}-X^j_{T_n})}{\gamma^2_\ep}\cr
=\int_0^1\lambda(dt)\int_0^1\lambda(ds)\Big(\bar k(t,s)-k(T_1,T_1)
\bar\alpha_1(t)\bar\alpha_1(s)\Big).\cr
}
$$
By using Theorem \ref{GD-G} one
gets the large deviation principle. Now, iterating the same
procedure up to $n$, one would achieve the following
(recall that $\Var(X^{j-1}_{T_j})=k_{j-1}(T_j,T_j)$):
$$
\displaylines{
\lim_{\ep\to 0} \E(\langle \lambda, X^n_{T_n+\ep \cdot}-X^n_{T_n}\rangle)
=0\cr
\lim_{\ep\to 0} \frac{\Var(\langle \lambda,X^n_{T_n+\ep
\cdot}-X^n_{T_n}\rangle)}{\gamma^2_\ep}
=\int_0^1\lambda(dt)\int_0^1\lambda(ds)\bar k_n(t,s)\cr
}
$$
with
$$
\bar k_n(t,s)=
\Big(\bar k(t,s)-
\sum_{j=1}^n k_{j-1}(T_j,T_j) \bar\alpha_j(t)\bar\alpha_j(s)\Big)
$$
Notice that, by (\ref{knbar}), $\bar k_n$ is
a continuous covariance function, being the (uniform) limit
of a continuous, symmetric and
positive definite function. Therefore, we can assert that
$\{(X^n_{T_n+\ep t}- X^n_{T_n})_{t\in[0,1]}\}_\ep$ satisfies a
large deviation principle on $C([0,1])$, with inverse speed
$\gamma^2_\ep$ and good rate function
$$
H_n(\varphi)=\left\{
\begin{array}{cl}
\displaystyle \frac 12\, \|\varphi\|^2_{\bar{\cl H}_n} &
\mbox{ if } \varphi\in \bar{\cl H}_n\\
+\infty & \mbox{ otherwise}
\end{array}
\right.
$$
Finally, since $X^n_{T_n+\ep t}=x_n+(X^n_{T_n+\ep t}-X^n_{T_n})$,
the large deviation principle as in the statement follows by contraction
and the associated  rate function is actually given by (\ref{Jn}).\cvd

\medskip

Before to continue with the asymptotic behavior of the $n$-fold
conditional bridge process, let us give some examples of
applications of last Theorem \ref{th-AB} to the fractional
Brownian motion and integrated Gaussian processes.

\subsection{Fractional Brownian motion}\label{FBM}

The following result holds as a consequence of Theorem
\ref{th-AB}:

\begin{theorem}\label{th-FBM}
Let $X$ be a fractional Brownian motion, with Hurst index
$H\in(0,1)$, and let $X^n$ denote the $n$-fold conditional process
as in (\ref{cond-process}). Then, the family of processes
$\{(X^n_{T_n+\ep t})_{t\in[0,1]}\}_\ep$ satisfies a
large deviation principle on $C([0,1])$, with inverse speed
$\ep^{2H}$ and good rate function
\begin{equation}\label{J-FBM}
J_n(h)=\left\{
\begin{array}{cl}
\displaystyle \frac 12\, \|h-x_n\|^2_{{\cl H}_H} & \mbox{ if }
h_0=x_n \mbox{ and } h-x_n\in {\cl H_H}\\
+\infty & \mbox{ otherwise}
\end{array}
\right.
\end{equation}
${\cl H}_H$ being the reproducing kernel Hilbert space associated to
the fractional Brownian motion itself.
\end{theorem}

Let us recall that a fractional Brownian motion $X$ with Hurst
index $H\in(0,1)$ is a continuous, non Markovian unless $H=1/2$, centered,
Gaussian process
whose covariance function is
$$
k_H(t,s)=\frac{t^{2H}+s^{2H}-|t-s|^{2H}}2.
$$

\sl Proof of Theorem \ref{th-FBM}\rm. We  show that both
Assumption \ref{ass-GD} and \ref{ass-A} do hold. First, one has
$$
\frac {\Cov(X_{T_n+\ep t}-X_{T_n}, X_{T_n+\ep
s}-X_{T_n})}{\ep^{2H}} =\Cov(X_t, X_s)
$$
because of the homogeneity and self-similarity properties
holding for the fractional Brownian motion,
so that the limit in (\ref{kbar}) trivially exists and
$\bar k(t,s)=k_H(t,s)$.
Concerning Assumption \ref{ass-A}, straightforward computations
(using Taylor expansion) allow easily to state that
$$
\lim_{\ep\to 0}\sup_{t\in[0,1]} \frac{|k_H(T_n+\ep t,T)-
k_H(T_n,T)|}{\ep^H} =0,
$$
for any $T>0$, so that $\bar\rho\equiv 0$. This in turn implies
that $\bar\alpha_j(t)=0$, for any $t\in[0,1]$ and $j=1,\ldots,n$,
as an immediate consequence of what developed in Lemma
\ref{bar} $(i)$. Then $\bar k_n\equiv k$ and
the statement now follows from Theorem \ref{th-AB}.
\cvd

Notice that the $n$-fold conditional fractional Brownian motion
satisfies a large deviation principle with the same rate function
as the non conditioned process. This means that the asymptotic
behavior of the $n$-fold conditional process does not depend on
the past, although $X$ is not Markovian unless $H=1/2$. This
is obvious for  $H=1/2$, when the process reduces to the
standard Brownian motion, and in some sense is quite natural when
$H<1/2$, that is when a short memory property holds.
Nevertheless, this is a little bit
surprising when $H>1/2$, a case in which the process has a long range
memory property.

\begin{example}\label{cheridito1}\rm
As an example, let us consider the process
$$
X_t=c B_t+c_H B^H_t,
$$
in which $c$ and $c_H$ are non null real numbers, $B$ stands for a standard
Brownian motion and $B^H$ denotes a fractional Brownian
motion with Hurst index $H\neq 1/2$. Suppose moreover that $B$ and $B^H$
are independent. Such a process has been
studied by Cheridito \cite{bib:cheridito1}, who proved
that $X$ is a semimartingale if and
only if $H\in (3/4,1)$, a property allowing to get
interesting applications in Finance. The covariance function
associated to $X$ is given by
$$
k(t,s)=c^2 k_{1/2}(t,s)+c_H^2 k_H(t,s).
$$
Then, by using arguments similar to
the ones developed in the proof of Theorem \ref{th-FBM},
one can state a large deviation principle for
$\{(X^n_{T_n+\ep t})_{t\in[0,1]}\}_\ep$ on $C([0,1])$, with inverse speed
$\ep^{2(H\wedge 1/2)}$ and good rate function associated to
the covariance function
\begin{equation}\label{sigma2H}
\bar k_n(t,s)=\sigma^2_H\, k_{H\wedge 1/2}(t,s),
\quad\mbox{with }
\sigma^2_H=\left\{
\begin{array}{cl}
c^2 & \mbox{ if } H>1/2\\
c_H^2 & \mbox{ if } H<1/2
\end{array}
\right.
\end{equation}
where $k_{H\wedge 1/2}$ denotes the covariance
function associated to a fractional Brownian motion with Hurst index
$H\wedge 1/2$. By contraction, the constant $\sigma^2_H$ can be
put inside the rate function, which becomes:
$$
J_n(h)=\left\{
\begin{array}{cl}
\displaystyle \frac 1{2\sigma^2_H}\, \|h-x_n\|^2_{{\cl H}_{H\wedge 1/2}} & \mbox{ if }
h_0=x_n \mbox{ and } h-x_n\in {\cl H}_{H\wedge 1/2}\\
+\infty & \mbox{ otherwise}
\end{array}
\right.
$$
${\cl H}_{H\wedge 1/2}$ being the reproducing kernel Hilbert space
associated to a fractional Brownian motion with Hurst index ${H\wedge 1/2}$.

\end{example}

\subsection{Integrated Gaussian Process}\label{IGP}

Let $Z$ be a centered  Gaussian process with covariance function
$\kappa (t,s)$ and let $X$ be the integrated process, i.e.,
\begin{equation}\label{X-IGP}
X_t=\int_0^t Z_u du.\end{equation} $X$ is a continuous, centered
Gaussian process whose covariance function $k$  is given by
\begin{equation}\label{k-IGP}
k(t,s)=\int_0^t \int_0^s \kappa(u,v) du dv.
\end{equation}

As a consequence of Theorem \ref{th-AB} one has:

\begin{theorem}\label{th-IGP}
Let $X$ be an integrated Gaussian process as in (\ref{IGP}), with
$\kappa(t,s)$ continuous, and let $X^n$ denote the $n$-fold
conditional process as in (\ref{cond-process}). Then, the family
$\{(X^n_{T_n+\ep t})_{t\in[0,1]}\}_\ep$
satisfies a large deviation principle on $C([0,1])$, with inverse
speed $\ep^{2}$ and good rate function
\begin{equation}\label{J-IGP}
J_n(h)=\left\{
\begin{array}{cl}
\displaystyle \frac 12\, \|h-x_n\|^2_{\bar{\cl H}_n} &
\mbox{ if } h_0=x_n \mbox{ and } h-x_n\in \bar{\cl H}_n\\
+\infty & \mbox{ otherwise}
\end{array}
\right.
\end{equation}
$\bar{\cl H}_n$ being the reproducing kernel Hilbert space
associated to the covariance function
$$
\bar k_n(t,s)=a_n^2\cdot ts, \quad\mbox{where}\quad
a_n^2=\kappa(T_n,T_n)-
\sum_{j=1}^n\frac{d_{j-1}(T_j)^2}{k_{j-1}(T_j,T_j)}
$$
and $d_{j-1}(T)$ is recursively defined as:
$d_0(T)=\int_0^T\kappa(T_n,u)du$ and as $i=1,2,\ldots, n-1$,
$$
d_i(T)=d_{i-1}(T)-\alpha_i(T)d_{i-1}(T_i)
$$
(recall that $k_j$ and $\alpha_j$ are defined through (\ref{kj})
and (\ref{alphaj}) respectively).

\end{theorem}

\proof Let us first show that Assumption \ref{ass-GD} holds,
with $\gamma_\ep=\ep$ and $\bar
k(t,s)=ts\,\kappa(T_n,T_n)$. In fact,
$$
\displaylines{ \Big|\frac 1{\ep^2}\,\Cov(X_{T_n+\ep t}-X_{T_n},
X_{T_n+\ep s}-X_{T_n})-ts\,\kappa(T_n,T_n)\Big|\leq\cr
\leq \frac
1{\ep^2}\,\int_{T_n}^{T_n+\ep}du
\int_{T_n}^{T_n+\ep}dv|\kappa(u,v)-\kappa(T_n,T_n)|\leq
\sup_{u,v\in[T_n, T_n+\ep]}|\kappa(u,v)-\kappa(T_n,T_n)| }
$$
and the last term goes to zero as $\ep\to 0$ because $\kappa$ is
continuous, then uniformly continuous on compact sets. Similarly,
one proves that also Assumption \ref{ass-A} holds, with
$\bar\rho(t,T)=t\,\int_0^T\kappa(T_n,v)dv$. The large
deviation principle is now an immediate application of Theorem
\ref{th-AB}. Finally, in order to give the above more explicit
expression for $\bar k_n$, we need the functions $\bar\alpha_j$.
By (\ref{alphajbar}), it is sufficient to show that
$$
\bar\rho_j(t,T)=d_j(T)\,t
$$
where $d_0(T)=\int_0^T\kappa(T_n,u)du$ and as $j=1,2,\ldots, n-1$,
$$
d_j(T)=d_{j-1}(T)-\alpha_j(T)d_{j-1}(T_i).
$$
We have already seen that $\bar\rho_0(t,T)=\bar\rho(t,T)=d_0(T)\,t$,
so that by (\ref{rhojbar}),
$$
\bar\rho_1(t,T)=\bar\rho_0(t,T)-\alpha_1(T)\bar\rho_0(t,T_1)=d_1(T)\,t,
$$
with $d_1(T)=d_{0}(T)-\alpha_1(T)d_{0}(T_1)$. By iteration, the
statement holds.
\cvd

\begin{remark}\label{cazzata}\rm
It follows the law of
an $n$-fold conditional integrated Gaussian process behaves asymptotically
as $a_nU\,t$, being $U$ a standard Gaussian random variable. Moreover,
a deeper view to the proof of Theorem \ref{th-IGP} shows that this
kind of ``degenerate'' behavior can be stated for any Gaussian process
whose
covariance function $k(t,s)$ is quite smooth, in particular if
both the first and the mixed second derivatives exist, the latter
being continuous on the diagonal points $(T,T)$. In fact, in this
case the asymptotic covariance $\bar k_n(t,s)$ for $X^n$ is again
of the type $a^2_n \cdot ts$.
\end{remark}

\begin{example}\label{mIBM-Xn}\rm \bf[$m$-fold integrated Brownian motion] \rm
Suppose that $X$ is defined as
$$
X_t=\int_0^tdu\Big(\int_0^{u}du_{m-1}\cdots\int_0^{u_2} du_1
W_{u_1}\Big),
$$
where $W$ denotes a standard Brownian motion.
It is known that $X$ is a centered, Gaussian process with
covariance function
$$
k(t,s)=\frac 1{(m!)^2}\int_0^{s\wedge t}(s-\xi)^m(t-\xi)^m\,d\xi
=\int_0^t\int_0^s \kappa(u,v)\,du\,dv,
$$
where
$$
\kappa(t,s)=\frac 1{((m-1)!)^2}\int_0^{t\wedge s}
(t-\xi)^{m-1}(s-\xi)^{m-1}\,d\xi
$$
(for details, see Chen and Li \cite{bib:cl}). Then, Theorem
\ref{IGP} applies to $X$. Notice that, for $T\leq T_n$ and $m\geq 1$,
$$
d_0(T)=\frac m{(m!)^2}\int_0^T(T_n-\xi)^{m-1}(T-\xi)^m\,d\xi.
$$

\end{example}

\begin{example}\label{IFBM-Xn}\rm \bf[Integrated  fractional Brownian motion] \rm
Suppose $X_t=\int_0^t Z_udu$, where $Z$ denotes a fractional Brownian
motion with Hurst index $H$. Then, the associated  covariance function
is
$$
k(t,s)=\int_0^t\int_0^s \kappa_H(u,v)\,du\,dv,\quad
\mbox{with } \kappa_H(t,s)=\frac 12 (t^{2H}+s^{2H} -|t-s|^{2H}).
$$
Again, Theorem \ref{IGP} immediately applies to $X$.
Here, for $T\leq T_n$, one has,
$$
d_0(T)=\frac 12 \Big[ T_n^{2H} T +\frac 1
{2H+1} (T^{2H+1}-T_n^{2H+1}-(T_n-T)^{2H+1})\Big] .
$$

\end{example}

\section{Large deviations for the bridge of the conditional
process}\label{bridge}

Let $(X^n_t)_{t\geq 0}$ be the $n$-fold conditional process defined
in Section \ref{model} and let us now consider the process  $Y^n$
defined as the {\it bridge} of the
process $X^n$, i.e, the process $X^n$ conditioned to be in
$y$ at the future time $T_n+\ep$. Then, in law one has,
\begin{equation}\label{bridge-process}
Y^n_{T_n+\ep t}=X^n_{T_n+\ep t}-\beta^\ep_{T_n+\ep t}(X^n_{T_n+\ep} -y),
\end{equation}
where
\begin{equation}\label{beta}
\beta^\ep_{T_n+\ep t}=\frac {k_n({T_n+\ep t},T_n+\ep)}{k_n(T_n+\ep,T_n+\ep)}.
\end{equation}
Now, in order to achieve a large deviation principle for
$\{(Y_{T_n+\ep t}^n)_{t\in[0,1]}\}_\ep$,
one needs a nice asymptotic behavior for $\beta^\ep_{T_n+\ep\cdot}$.
In fact, one has

\begin{lemma}\label{lemma-beta}
Let Assumptions \ref{ass-GD} and \ref{ass-A} be satisfied. Then
there exists the limit
\begin{equation}\label{betabar}
\bar\beta_t=\lim_{\ep\to 0} \beta^\ep_{T_n+\ep t}
=\frac{\bar k_n(t,1)}{\bar k_n(1,1)},\quad
\mbox{uniformly as } t\in[0,1].
\end{equation}
\end{lemma}
\proof One has,
$$
\displaylines{  |\beta^\ep_{T_n+\ep t}-\bar\beta_t|=
 \Big|\frac{k_n(T_n+\ep t,T_n+\ep)}{k_n(T_n+\ep ,T_n+\ep
 )}-\frac{\bar k_n(t,1)}{\bar k_n(1,1)}\Big|\cr\leq\!\!
\frac {\gamma^2_\ep}{k_n(T_n+\ep  ,T_n+\ep )}\Big| \frac
{k_n(T_n+\ep t,T_n+\ep )}{\gamma^2_\ep} -\bar k_n(t,1)\Big|+
|\bar k(t,1)| \Big|\frac {\gamma^2_\ep}{k_n(T_n+\ep ,T_n+\ep
)}-\frac 1{\bar k_n(1,1)}\Big|}
$$
From (\ref{knbar}),
$$
\displaylines{
\lim_{\ep \to 0} \sup_{t\in[0,1]}\Big| \frac
{k_n(T_n+\ep t,T_n+\ep )}{\gamma^2_\ep} -\bar k_n(t,1)\Big|=0
\mbox{ and }
\lim_{\ep \to 0}\Big| \frac{\gamma^2_\ep}
{k_n(T_n+\ep ,T_n+\ep )} -\frac 1{\bar k_n(1,1)}\Big|=0,\cr
}
$$
so that the statement holds.     \cvd

It is now easy to prove a first large deviation principle. But,
as we will see, there are cases in which the next, immediate,
result turns out to be degenerate in some sense.
So, let us split this
Section in two part, the former  containing a first result and the
latter developing some refinements.

\subsection{A first large deviation result for the bridge}\label{immediate-LD}

\begin{theorem}\label{th-bridge}
Let $Y^n$ be the bridge of the $n$-fold conditional process $X^n$,
as defined in (\ref{bridge-process}). Under
Assumptions \ref{ass-GD} and
\ref{ass-A},  the family of processes $\{(Y^n_{T_n+\ep
t})_{t\in[0,1]}\}_\ep$ satisfies a large deviation
principle on $C([0,1])$, with inverse speed $\gamma^2_\ep$ and
good rate function
\begin{equation}\label{J-bridge-0}
J_Y(h)=\left\{
\begin{array}{cl}
\displaystyle \frac 12\, \|h-\bar m\|^2_{\bar{\cl H}_Y}
& \mbox{ if } h_0=x_n, h_1=y, h-\bar m\in \bar{\cl H}_Y\\
+\infty & \mbox{ otherwise}
\end{array}
\right.
\end{equation}
where $\bar m_t =x_n+\bar\beta_t(y-x_n)$ and
$\bar{\cl H}_Y$ is the reproducing kernel Hilbert space
associated to the covariance function
$$
\bar k_Y(t,s)=\bar k_n(t,s)-\bar\beta_s\bar k_n(t,1) =\bar
k_n(t,s)-\frac{\bar k_n(t,1)\bar k_n(s,1)}{\bar k_n(1,1)}.
$$

\end{theorem}

\proof
First, let us set
$$
U^n_{T_n+\ep t}=Y^n_{T_n+\ep t}-\bar m_t,
\quad\mbox{where }
\bar m_t
=x_n+\bar\beta_t(y-x_n)
=\lim_{\ep\to 0}\E(Y^n_{T_n+\ep t})
$$
and notice that, by (\ref{bridge-process}), (\ref{mnbar}) and
(\ref{betabar}), the above limit holds uniformly as $t\in[0,1]$.
We will start by showing a large deviation principle for
$\{U^n_{T_n+\ep \cdot}\}_\ep$, by using again Theorem \ref{GD-G}. In
fact,
$$
\lim_{\ep\to 0}\E\Big(\<\lambda,U^n_{T_n+\ep \cdot}\>\Big)=\int_0^1
\lambda(dt)\E(U^n_{T_n+\ep t})=0,
$$
for any $\lambda\in\cl M[0,1]$. Moreover,
from (\ref{knbar}) and Lemma \ref{lemma-beta}, one has
$$
\displaylines{
\lim_{\ep\to 0}
\frac{\Cov(U^n_{T_n+\ep t}, U^n_{T_n+\ep s})}{\gamma^2_\ep}
=\lim_{\ep\to 0}
\frac{\Cov(Y^n_{T_n+\ep t}, Y^n_{T_n+\ep s})}{\gamma^2_\ep}=\cr
=\lim_{\ep\to 0}
\frac{\Cov(X_{T_n+\ep t}^n,
X_{T_n+\ep s}^n)-\beta^\ep_{T_n+\ep s}\Cov(X_{T_n+\ep t}^n,
X_{T_n+\ep}^n)}{\gamma^2_\ep}=\cr
=\bar k_n(t,s)-\bar\beta_s\bar k_n(t,1)
=:\bar k_Y(t,s),\cr
}
$$
uniformly as $s,t\in[0,1]$, so that
$$
\displaylines{
\lim_{\ep\to 0}
\frac{
\Var\Big(\<\lambda, U^n_{T_n+\ep
\cdot}\>\Big)}{\gamma^2_\ep}\cr
=\displaystyle\lim_{\ep\to 0}
\frac{1}{\gamma^2_\ep}\int_0^1\lambda(dt)\int_0^1\lambda(ds)\,
\Cov(U^n_{T_n+\ep t} U^n_{T_n+\ep s})
=\int_0^1\lambda(dt)\int_0^1\lambda(ds)\,\bar k_Y(t,s),\cr
}
$$
for any $\lambda\in\cl M[0,1]$.
We can then assert that the family of processes $\{(U^n_{T_n+\ep
t})_{t\in[0,1]}\}_\ep$ does satisfy a large deviation
principle on $C([0,1])$, with inverse speed $\gamma^2_\ep$ and
good rate function
$$
J_U(\varphi)=\left\{
\begin{array}{cl}
\displaystyle \frac 12\, \|\varphi\|^2_{\bar{\cl H}_Y}
& \mbox{ if } \varphi_0=\varphi_1=0, \varphi\in \bar{\cl H}_Y\\
+\infty & \mbox{ otherwise}
\end{array}
\right.
$$
being $\bar{\cl H}_Y$ the reproducing kernel Hilbert space associated to the
covariance function
$$
\bar k_Y(t,s)=\bar k_n(t,s)-\bar\beta_s\bar k_n(t,1)
=\bar k_n(t,s)-\frac{\bar k_n(t,1)\bar k_n(s,1)}{\bar k_n(1,1)}.
$$
Let us stress that the condition $\varphi_0=\varphi_1=0$ is
trivially satisfied if $\varphi\in \bar{\cl H}_Y$ (it immediately
follows from the fact that $k_Y(0,s)=k_Y(1,s)=0$ for any $s$), but
we have chosen to write it for the sake of clearness. Now, since
$Y^n_{T_n+\ep t}=U^n_{T_n+\ep t}+\bar m_t$, by contraction one
immediately obtains the large deviation principle for
$\{(Y^n_{T_n+\ep t})_{t\in[0,1]}\}_\ep$ on $C([0,1])$, with
inverse speed $\gamma^2_\ep$ and good rate function as in
(\ref{J-bridge-0}). \cvd

\begin{remark}\rm
The rate function $J_Y$ given by (\ref{J-bridge-0})
can also be written in the
following way:
\begin{equation}\label{J-bridge-1}
J_Y(h)=\left\{
\begin{array}{cl}
\displaystyle \frac 12\,\Big( \|h-x_n\|^2_{\bar{\cl
H}_n}-\frac{(y-x_n)^2}{\bar k_n(1,1)}\Big)
& \mbox{ if } h_0=x_n, h_1=y, h-x_n\in \bar{\cl H}_n\\
+\infty & \mbox{ otherwise}
\end{array}
\right.
\end{equation}
$\bar{\cl H}_n$ being the reproducing kernel Hilbert space
associated to the covariance function $k_n$ defined in
(\ref{kbar-n}).  Such a representation agrees with well known formulas,
 for example whenever $X$ is a standard
Brownian motion (see e.g. Baldi, Caramellino and Iovino \cite{bib:bci}).
The proof of (\ref{J-bridge-1}) is postponed to Appendix A.

\end{remark}

\begin{example}\label{FBM-Y}\rm \bf [Fractional
Brownian motion] \rm  Following Section \ref{FBM},
let $X$ be a fractional
Brownian motion with Hurst index $H$ and let
$X^n$ be the associated  $n$-fold conditional process.
As seen in Theorem \ref{th-FBM}, both Assumption \ref{ass-GD} and \ref{ass-A}
hold and the asymptotic covariance
function $\bar k_n(t,s)$ coincides with the original one $k_H(t,s)$.
By applying Theorem \ref{th-bridge}, the bridge process $Y^n$
satisfies a functional large deviation principle
for small time, with inverse speed $\ep^{2H}$ and
good rate function
\begin{equation}\label{ponte-MBF0}
J_Y(h)=\left\{
\begin{array}{cl}
\displaystyle \frac 12 \|h-\bar m\|^2_{\bar{\cl H}_Y} & \mbox{
if } h_0=x_n, h_1=y,
h-x_n\in \bar{\cl H}_Y\\
+\infty & \mbox{ otherwise}
\end{array}
\right.
\end{equation}
$\bar{\cl H}_Y$ being the reproducing kernel Hilbert space
associated to the covariance function
$$
\bar k_Y(t,s)=k_H(t,s)-k_H(t,1)k_H(1,s).
$$
By using (\ref{J-bridge-1}), $J_Y$ can be written also in terms of
the reproducing kernel Hilbert space ${\cl H}_H$
associated to the original fractional
Brownian motion $X$:
\begin{equation}\label{ponte-MBF1}
J_Y(h)=\left\{
\begin{array}{cl}
\displaystyle \frac 12\Big(
\|h-x_n\|^2_{{\cl H}_H} -(y-x_n)^2\Big)& \mbox{ if } h_0=x_n, h_1=y,
h-x_n\in {\cl H}_H\\
+\infty & \mbox{ otherwise}
\end{array}
\right.
\end{equation}
Whenever $H=1/2$, that is $X$ is a standard
Brownian motion, then the above result is well known and widely applied
in the literature.
Let us stress that, as well as the $n$-fold conditional fractional
Brownian motion, also its bridge satisfies a large deviation
principle which is independent of all the past except for what
happens at time $T_n$.
\end{example}

\medskip

\begin{example}\rm \label{cheridito2}
{\bf [Cheridito process]} Let $X$ be the process as in Example
\ref{cheridito1}: $X_t=c B_t+c_H B^H_t$, with $c, c_H\neq 0$
constant numbers, $B$ and $B^H$ are independent,
with $B$ a Brownian motion and $B^H$ a fractional
Brownian motion with Hurst index $H\neq 1/2$.
By developing arguments similar
to the ones in Example \ref{FBM-Y}, one obtains that
the bridge of the associated $n$-fold conditional
process satisfies a large deviation principle.
By taking into account the results
in Example \ref{cheridito1} and formula (\ref{J-bridge-1}), one easily obtains that
the inverse speed is equal to
$\ep^{2(H\wedge 1/2)}$ and the good rate function is given
by the following formula:
\begin{equation}\label{ponte-cheridito}
J_Y(h)=\left\{
\begin{array}{cl}
\displaystyle \frac 1{2\sigma^2_H}\Big(
\|h-x_n\|^2_{{\cl H}_{H\wedge 1/2}} -(y-x_n)^2\Big)& \mbox{ if } h_0=x_n, h_1=y,
h-x_n\in {\cl H}_{H\wedge 1/2}\\
+\infty & \mbox{ otherwise}
\end{array}
\right.
\end{equation}
where $\sigma^2_H$ is given by (\ref{sigma2H}).

\end{example}

\medskip

\begin{example} \rm \label{IGP-Y-deg}\bf [Integrated Gaussian processes]
\rm  Following Section \ref{IGP},
let $X^n$ be the $n$-fold conditional process
when $X$ is an integrated Gaussian process as in (\ref{X-IGP}).
Under the hypotheses of Theorem \ref{th-IGP},
Assumptions \ref{ass-GD} and
\ref{ass-A} hold, and a  functional large deviation principle
for $X^n$ follows, with asymptotic covariance function
$\bar k_n(t,s)=a_n^2\cdot ts$ for a
suitable constant $a_n^2$. Now, by applying Theorem \ref{th-bridge},
one obtains a functional large deviation principle for the bridge
process $Y^n$ as well, but unfortunately one gets
a degenerate asymptotic behavior because the associated rate function
turns out to be
$$
J_Y(h)=\left\{
\begin{array}{cl}
\displaystyle 0& \mbox{ if } h=\bar m\\
+\infty & \mbox{ otherwise}.
\end{array}
\right.
$$
This follows from the fact that, since $\bar{\cl H}_n$ is
``spanned'' by the covariance function $\bar k_n(t,s)=a^2_n\cdot
ts$, it contains only the paths running at constant speed.
Then, $J_Y$ is finite only for $h$
such that $h-\bar m=ct$. Since here  $\bar m_t=x_n+(y-x_n)t$, the additional constraints $h_0=x_n$ and
$h_1=y$ give the unique path $h=\bar m$.

\end{example}

Notice that Theorem \ref{th-bridge} gives
an unsatisfactory large deviation result not only for integrated Gaussian
processes but also for Gaussian processes whose (original)
covariance function
is smooth enough: as observed in Remark \ref{cazzata}, in this
case the asymptotic covariance function is $const\cdot ts$ as well, and the
same degenerate behavior holds for the rate function.
This motivates next section,
in which we study some refinements allowing one to state non trivial
large deviation estimates, or more precisely the right large deviation speed.

\subsection{Faster large deviations
for the bridge}\label{refined-LD}

In this Section we
prove a refined version of Theorem \ref{th-bridge}: we study here
the exact (faster) speed giving a non trivial rate function
whenever the covariance is smooth.
With the same notation of Section \ref{model},
Assumptions \ref{ass-GD} and \ref{ass-A} must be
here strengthened as follows:
\begin{assumption}\label{ass-ref1}
For some $\alpha\in(0,1]$,
\begin{itemize}
\item[$(i)$]
there exist a function $\bar\varphi(t,s)$, a
constant $a^2$ and a remaining term  $\cl R^1_\ep(t,s)$ (depending on
$T_n$) such that
\begin{equation}\label{kbarref1}
\begin{array}{l}
\Cov(X_{T_n+\ep t}-X_{T_n}, X_{T_n+\ep s}-X_{T_n})=\ep^2\left[a^2\,
ts + \bar\varphi(t,s)\ep^{\alpha}+\cl R^1_\ep(t,s)\right]\smallskip\\
\mbox{ with } \displaystyle\lim_{\ep\to 0} \sup_{s,t\in[0,1]} \frac{|\cl
R^1_\ep(t,s)|}{\ep^{\alpha}}=0;
\end{array}
\end{equation}
\item[$(ii)$]
for any fixed $T>0$, there exist a function $\bar\psi(t,T)$,
a constant $c(T)$ and a remaining term  $\cl R^2_\ep(t;T)$ (depending on
$T_n$) such that
\begin{equation}\label{rhobarref1}
\begin{array}{l}
k(T_n+\ep t,T)-k(T_n,T)=\ep\left[c(T)\, t+\bar \psi
(t;T)\ep^{\alpha}+\cl R^2_\ep(t;T)\right]\smallskip\\
\mbox{ with }\displaystyle \lim_{\ep\to 0} \sup_{t\in[0,1]} \frac{|\cl
R^2_\ep(t;T)|}{\ep^{\alpha}}=0.
\end{array}
\end{equation}
\end{itemize}

\end{assumption}

As a consequence of Assumption \ref{ass-ref1},
by using the same arguments as in Lemma
\ref{bar} one immediately prove the following
\begin{lemma}\label{lemma-refined1}
For $j=1,\dots, n$
\begin{equation}
k_j(T_n+\ep t)-k_j(T_n,T)=\ep[c_j(T)\, t +\bar \psi_j
(t;T)\ep^{\alpha}+\cl R^2_\ep(t;T)],
\end{equation}
where setting, $c_0\equiv c$, and $\bar\psi_0\equiv \bar\psi$,
$c_j$ and $\bar\psi_j$ are given by
$$
c_j(T)=c_{j-1}(T)-\alpha_j(T) c_{j-1}(T_j)\quad\mbox{and}\quad
\bar\psi_j(t;T)=\bar\psi_{j-1}(t;T)-\alpha_j(T)\bar\psi_{j-1}(t;T_j).
$$
 Moreover
$$\Cov(X^j_{T_n+\ep t}-X^j_{T_n},X^j_{T_n+\ep s}-X^j_{T_n})=
\ep^2[ a^2_j\, ts +\bar \varphi_j(t,s)\ep^{\alpha}+\cl
R^{1,j}_\ep(t,s)],
$$
where ${\cl R}^{1,j}_\ep(t,s)\to 0$ as $\ep\to 0$
uniformly on $[0,1]\times[0,1]$,  $a_j=a-\sum_{\ell=1}^j c_{\ell-1}^2(T_\ell)$
and
\begin{equation}\label{barphij1}
\bar\varphi_j(t,s)=\bar\varphi(t,s)-\sum_{\ell=1}^j
\frac{c_{\ell-1}(T_\ell)}
{k_{\ell-1}(T_\ell,T_\ell)}(\bar\psi_{\ell-1}(t;T_\ell) \, s
+\bar\psi_{\ell-1}(s;T_\ell) \, t).
\end{equation}

In particular, since $X^n_{T_n}=x_n$,
\begin{equation}\label{barknref1}
k_n(T_n+\ep t,T_n+\ep s)=\ep^2[ a^2_n\, ts +\bar
\varphi_n(t,s)\ep^{\alpha}+\cl R^{1,n}_\ep(t,s)],
\mbox{ with }
\lim_{\ep\to0}\sup_{t,s\in[0,1]}\frac{|\cl R^{1,n}_\ep(t,s)| }{\ep^\alpha}=0.\end{equation}
\end{lemma}

Then, one has

\begin{theorem}\label{th-bridgeref1}
Let $Y^n$ be the bridge of the $n$-fold conditional process $X^n$, as
defined in (\ref{bridge-process}). If Assumption \ref{ass-ref1}
holds, then  the family of processes $\{(Y^n_{T_n+\ep
t})_{t\in[0,1]}\}_\ep$ satisfies a large deviation
principle on $C([0,1])$, with inverse speed $\ep^{2+\alpha}$ and
good rate function
\begin{equation}
 J_Y(h)=\left\{
\begin{array}{cl}
\displaystyle \frac 12\, \|h-\bar m\|^2_{\bar{\cl  H}_Y} & \mbox{ if } h-
\bar m\in \bar{\cl H}_Y\\
+\infty & \mbox{ otherwise}
\end{array}
\right.
\end{equation}
where $\bar m_t=x_n+\bar\beta_t(y-x_n)$ and
$\bar{\cl  H}_Y$ is the reproducing kernel Hilbert space
associated to the covariance function
\begin{equation}\label{barkYref1}
\bar k_Y(t,s)=\bar \varphi_n(t,s)+ t s\, \bar \varphi_n(1,1)- t\,
\bar \varphi_n(1,s)-s \,\bar \varphi_n(t,1).
\end{equation}
\end{theorem}
\proof The proof is the same as  Theorem \ref{th-bridge}. It is
enough to observe that in this case one has, from (\ref{barknref1}),
$$
\displaylines{\Cov(Y^n_{T_n+\ep t}, Y^n_{T_n+\ep s})=k_n(T_n+\ep
t,T_n+\ep s)- \frac{k_n(T_n+\ep t, T_n+\ep)k_n(T_n+\ep , T_n+\ep
s)}{k_n(T_n+\ep , T_n+\ep)}\cr
=\ep^2\Big[ a^2_n\, ts +\bar
\varphi_n(t,s)\ep^{\alpha}-\frac{( a^2_n\, t +\bar
\varphi_n(t,1)\ep^{\alpha})( a^2_n\, s +\bar
\varphi_n(1,s)\ep^{\alpha})}{ a^2_n\,  +\bar
\varphi_n(1,1)\ep^{\alpha}+\cl R^{1,n}_\ep(1,1)}+ \cl
R^{1,n}_\ep(t,s)\Big]\cr
=\!\frac{a^2_n\ep^2}{ a^2_n\,  +\bar
\varphi_n(1,1)\ep^{\alpha}+\cl R^{1,n}_\ep(1,1)}\Big(\!(\bar
\varphi_n(t,s)+ t s\, \bar \varphi_n(1,1)- t\, \bar
\varphi_n(1,s)-s \,\bar \varphi_n(t,1))\ep^\alpha +\cl
R^{1,n}_\ep(t,s)\!\Big).
}
$$
Therefore
\begin{equation}\label{app-cov}
\lim_{\ep\to 0} \frac{\Cov(Y^n_{T_n+\ep t}, Y^n_{T_n+\ep s})}{\ep^{2+\alpha}}=
\bar \varphi_n(t,s)+ t s\, \bar \varphi_n(1,1)- t\, \bar
\varphi_n(1,s)-s \,\bar \varphi_n(t,1),
\end{equation}
uniformly as $s,t\in[0,1]$ and the thesis holds. \cvd

\begin{remark}\label{false-bridge}\rm
Notice that $\bar\varphi_n$ is symmetric and continuous whereas
it is not
positive definite in general, so that
it is not necessarily a covariance function. Nevertheless, $k_Y$ given by
(\ref{barkYref1}) does represent a covariance function, as an
immediate consequence of (\ref{app-cov}). However, if
$\bar\varphi_n$ was a covariance function, a curious effect would happen:
the asymptotic behavior of the bridge is regulated by
a covariance function which coincides with the one
associated to what is usually called ``the false bridge'',
that is a process
of the type $Z_t-tZ_1$, being $Z$ a Gaussian
process with covariance $\bar\varphi_n$.

\end{remark}

Now, if the function $k(t,s)$ is more regular, then Theorem
\ref{th-bridgeref1} would give again a degenerate behavior.
In fact, suppose that  $k$ has continuous derivatives up to
the third order. Then, since $k$ is symmetric, by straightforward
computations  one obtains
$$
\bar\varphi(t,s)=\frac 1{3!}(3\partial^3_{tts} k(T_n,T_n) t^2 s+
3\partial^3_{tss} k(T_n,T_n) t s^2)=\frac 12 \partial^3_{tts}
k(T_n,T_n)ts( t+ s)
$$
and by using (\ref{barphij1}) one arrives to show that
$\bar\varphi_n(t,s)=b_n ts(t +s)$, for a suitable constant $b_n$.
Therefore, by (\ref{barkYref1}) one has
$$
\bar k_Y(t,s)=b_nts(t +s)+2b_n ts- b_n(t^2 +t)s- b_n(s+s^2)t\equiv 0,
$$
and again a trivial large deviation principle holds
for the bridge of the conditional process. Let us refine further on the
hypothesis.

\begin{assumption}\label{ass-ref2}
For some  $\alpha\in(0,1]$,
\begin{itemize}
\item[$(i)$]
there exist a function $\bar\varphi(t,s)$, constants $a^2$ and $b$
and a remaining term $\cl R^1_\ep(t,s)$
(depending on $T_n$) such that
\begin{equation}\label{kbarref2}
\begin{array}{l}
\Cov(X_{T_n+\ep t}-X_{T_n}, X_{T_n+\ep s}-X_{T_n})=\smallskip\\
\quad\quad=\ep^2[a^2\, ts + b (t^2s+ts^2) \ep+
\bar\varphi(t,s)\ep^{1+\alpha}
 +\cl R^1_\ep(t,s)]\smallskip\\
\quad\quad\mbox{ with }\lim_{\ep\to 0} \sup_{s,t\in[0,1]} \frac{|\cl
R^1_\ep(t,s)|}{\ep^{1+\alpha}}=0;
\end{array}
\end{equation}
\item[$(ii)$]
for any fixed $T>0$, there exist a  function $\bar\psi(t,T)$, constants
$c(T)$ and $d(T)$ and a remaining term $\cl R^2_\ep(t;T)$
 (depending on $T_n$) such that
\begin{equation}\label{rhobarref2}
\begin{array}{l}
k(T_n+\ep t,T)-k(T_n,T)=\ep[c(T)\, t+ d(T) \, t^2 \ep+\bar \psi
(t;T) \ep^{1+\alpha}+\cl R^2_\ep(t;T)]\smallskip\\
\mbox{ with }\displaystyle \lim_{\ep\to 0} \sup_{t\in[0,1]} \frac{|\cl
R^2_\ep(t;T)|}{\ep^{1+\alpha}}=0.
\end{array}
\end{equation}
\end{itemize}

\end{assumption}

Let us remark that if $k(t,s)$ is smooth enough then
one immediately has  $a^2=\partial^2_{ts}k(T_n,T_n)$, $b=\frac 12
\partial^3_{tts} k(T_n,T_n)$, $c(T)=\partial_t k(T_n,T)$ and
$d(T)=\frac 12 \partial^2_{tt}k(T_n,T)$.

Moreover, as an immediate consequence of Assumption \ref{ass-ref2},
by using the same arguments as in Lemma
\ref{bar} and Lemma \ref{lemma-refined1}, one
proves that
\begin{lemma}\label{lemma-refined2}
For $j=1,\dots, n$
\begin{equation}
k_j(T_n+\ep t)-k_j(T_n,T)=\ep[c_j(T)\, t + d_j(T) \, t^2 \ep +\bar
\psi_j (t;T)\ep^{1+\alpha}+\cl R^{2,j}_\ep(t;T)]
\end{equation}
where setting, $c_0(T)=c(T)$, $d_0(T)=d(T)$  and $\bar\psi_0\equiv
\bar\psi$, $c_j(T)$,  $d_j(T)$ and $\bar\psi_j$ are defined in the
following way,
$$
\displaylines{
c_j(T)=c_{j-1}(T)-\alpha_j(T) c_{j-1}(T_j), \quad
d_j(T)=d_{j-1}(T)-\alpha_j(T) d_{j-1}(T_j),\cr
\bar\psi_j(t;T)=\bar\psi_{j-1}(t;T)-\alpha_j(T)\bar\psi_{j-1}(t;T_j).\cr
}
$$
 Moreover
$$\Cov(X^j_{T_n+\ep t}-X^j_{T_n},X^j_{T_n+\ep s}-X^j_{T_n})=
\ep^2[ a^2_j\, ts + b_j(t^2 s +t s^2)\ep+\bar
\varphi_j(t,s)\ep^{1+\alpha}+\cl R^{1,j}_\ep(t,s)],
$$
where  $a^2_j=a^2-\sum_{\ell=1}^j c_{\ell-1}^2(T_\ell)$,
$b_j=b-\sum_{\ell=1}^j c_{\ell-1}(T_\ell)d_{\ell-1}(T_\ell)$ and
\begin{equation}\label{barphij2}
\bar\varphi_j(t,s)=\left\{
\begin{array}{ll}
\displaystyle\bar\varphi(t,s)-\sum_{\ell=1}^j
\frac{c_{\ell-1}(T_\ell)}
{k_{\ell-1}(T_\ell,T_\ell)}(\bar\psi_{\ell-1}(t;T_\ell) \, s
+\bar\psi_{\ell-1}(s;T_\ell) \, t) & \mbox{for}\, \alpha<1\\
\displaystyle\bar\varphi(t,s)-\sum_{\ell=1}^j
\frac{c_{\ell-1}(T_\ell)}
{k_{\ell-1}(T_\ell,T_\ell)}(\bar\psi_{\ell-1}(t;T_\ell) \, s
+\bar\psi_{\ell-1}(s;T_\ell) \, t)+\\
\qquad\quad-\displaystyle\sum_{\ell=1}^j
\frac{d_{\ell-1}^2(T_\ell)} {k_{\ell-1}(T_\ell,T_\ell)}t^2 s^2 &
\mbox{for}\, \alpha=1.
\end{array}
\right.
\end{equation}

In particular since $X^n_{T_n}=x_n$, one has
\begin{equation}\label{barknref2}k_n(T_n+\ep t,T_n+\ep s)=\ep^2[ a_n\, ts + b_n (t^2 s+t
s^2)\ep +\bar \varphi_n(t,s)\ep^{1+\alpha}+\cl
R^{1,n}_\ep(t,s)].
\end{equation}

\end{lemma}

Let us stress that in Lemma \ref{lemma-refined2}, the notation
${\cl R}_\ep$ (with some suitable superscript)
stands for a generical remaining term, which
uniformly converges to $0$ as $\ep\to 0$.

Then, one has
\begin{theorem}\label{th-bridgeref2}
Let $Y^n$ be the bridge of the $n$-fold conditional process $X^n$,
as defined in (\ref{bridge-process}). If Assumption \ref{ass-ref2}
holds, then  the family of processes $\{(Y^n_{T_n+\ep
t})_{t\in[0,1]}\}_\ep$ satisfies a large deviation
principle on $C([0,1])$, with inverse speed $\ep^{3+\alpha}$ and
good rate function
\begin{equation}
 J_Y(h)=\left\{
\begin{array}{cl}
\displaystyle \frac 12\, \|h-\bar m\|^2_{\bar{\cl  H}_Y} & \mbox{ if } h-\bar m\in \bar{\cl H}_Y\\
+\infty & \mbox{ otherwise}
\end{array}
\right.
\end{equation}
where $\bar m_t=x_n+\bar\beta_t(y-x_n)$ and
$\bar{\cl  H}_Y$ is the reproducing kernel Hilbert space
associated to the covariance function
\begin{equation}\label{barkyref2}
\bar k_Y(t,s)=\left\{
\begin{array}{ll}
\bar \varphi_n(t,s)+ t s\, \bar \varphi_n(1,1)- t\,
\bar \varphi_n(1,s)-s \,\bar \varphi_n(t,1)&\mbox{if }\alpha<1\smallskip\\
\begin{array}{l}
b_n^2( t s^2+ t^2s-t^2 s^2 -st)+  \smallskip\\
+\bar\varphi_n(t,s)+ t s\, \bar \varphi_n(1,1)- t\, \bar
\varphi_n(1,s)-s \,\bar \varphi_n(t,1)
\end{array}
&\mbox{if }\alpha=1
\end{array}
\right.
\end{equation}

\end{theorem}

\proof
The proof is the same as  Theorem \ref{th-bridge}.
It is enough to observe that
$$
\displaylines{\Cov(Y^n_{T_n+\ep t}, Y^n_{T_n+\ep s})=k_n(T_n+\ep
t,T_n+\ep s)- \frac{k_n(T_n+\ep t, T_n+\ep)k_n(T_n+\ep , T_n+\ep
s)}{k_n(T_n+\ep , T_n+\ep)}=\cr \ep^2\Big[ a^2_n\, ts + b_n (t^2 s+t
s^2)\ep +\bar \varphi_n(t,s)\ep^{1+\alpha}\cr-\frac{( a^2_n\, t +
b_n (t^2 +t )\ep +\bar \varphi_n(t,1)\ep^{1+\alpha})( a^2_n\, s +
b_n ( s+ s^2)\ep +\bar \varphi_n(1,s)\ep^{1+\alpha})}{ a^2_n\,  +
2b_n \ep +\bar \varphi_n(1,1)\ep^{1+\alpha}+\cl R^1_\ep(1,1)}+ \cl
R^1_\ep(t,s)\Big]=\cr \frac{a^2_n\ep^2}{  a^2_n\,  + 2b_n \ep +\bar
\varphi_n(1,1)\ep^{1+\alpha}+\cl R^1_\ep(1,1)}\Big((\bar
\varphi_n(t,s)+ t s\, \bar \varphi_n(1,1)- t\, \bar
\varphi_n(1,s)\cr
-s \,\bar \varphi_n(t,1))\ep^{1+\alpha }+b_n^2(
t s^2+ t^2s-t^2 s^2 -st)\ep^2+\cl R^1_\ep(t,s)\Big),}
$$ therefore
 the thesis holds.  \cvd

Let us observe that if the covariance function $k(t,s)$ is more
regular, that is  $C^{4+\beta}$ for some $\beta\geq 0$, then
 Theorem \ref{th-bridgeref2} continues to hold
 and  the associated asymptotic covariance $\bar k_Y$, given
by (\ref{barkyref2}), is not in general degenerate.  In
fact, since  the fourth derivatives exist, one gets
$\bar\varphi_n(t,s)=e_n (t^3 s+ t s^3) +f_n t^2 s^2$. Therefore,
tedious but straightforward computations will give $\bar k_Y(t,s)
=const  \,\, ts (1-t)(1-s)$.

\bigskip

Let us now come back to the example suggesting to refine our first result
for the bridge of the $n$-fold conditional process, that is the
integrated Gaussian
process:
$$
X_t=\int_0^t Z_u du,
$$
$Z$ being a centered Gaussian process with covariance function
$\kappa(t,s)$. We are looking for conditions on
$\kappa$ so that Assumption
\ref{ass-ref1} or \ref{ass-ref2} is satisfied
and then a large deviation
principle as in  Theorem \ref{th-bridgeref1} or
\ref{th-bridgeref2} does hold.
One has

\begin{proposition}\label{prop-IGPref12}

1. Suppose that for some $\alpha\in(0,1]$,
$$
\begin{array}{rl}
\kappa(T_n+\ep u, T_n+\ep v)&=\kappa(T_n,T_n) +\ep^\alpha \hat
g(u,v) + \hat{\cl R}_\ep(u,v)\smallskip\\
\displaystyle\int_0^T \kappa(T_n+\ep u, v)dv&
\displaystyle= \int_0^T dv \,  \kappa(T_n,v)
+\ep^\alpha \tilde  g(u;T) + \tilde{\cl R}_\ep(u;T),\quad T>0,
\end{array}
$$
(the above functions and remaining terms may all depend
on $T_n$) with  $\hat g\in L^1([0,1]^2)$, $\tilde g(\cdot;T)\in L^1([0,1])$
and
$$
\lim_{\ep\to 0} \ep^{-\alpha}\,\|\hat{\cl
R_\ep}(\cdot,\cdot)\|_{L^1([0,1]^2)}= 0,\quad
\lim_{\ep\to 0} \ep^{-\alpha}\,\|\tilde{\cl
R_\ep}(\cdot;T)\|_{L^1([0,1])}= 0.
$$
Then, Assumption \ref{ass-ref1} holds, with
$$
\bar \varphi(t,s)=\int_0^t du \int_0^s dv \, \hat
g(u,v)\quad\mbox{and}\quad \bar \psi(t,T)=\int_0^t du \, \tilde
g(u;T).
$$

2. Suppose that for some $\alpha\in(0,1]$,
$$
\begin{array}{rl}
\kappa(T_n+\ep u, T_n+\ep v)&=\kappa(T_n,T_n) + \ep
e\cdot(u+v)+\ep^{1+\alpha}\hat g(u,v) + \hat{\cl R}_\ep(u,v)\smallskip\\
\displaystyle\int_0^T \kappa(T_n+\ep u, v)dv&=\displaystyle \int_0^T dv \,  \kappa(T_n,v)
+\ep u f(T) +\ep^{1+\alpha }\tilde g(u;T) + \tilde{\cl
R}_\ep(u;T), T>0,
\end{array}
$$
(the above functions, remaining terms and the constants $e$ and $f(T)$ may all depend
on $T_n$), with  $\hat g\in L^1([0,1]^2)$, $\tilde g(\cdot;T)\in L^1([0,1])$
and
$$
\lim_{\ep\to 0} \ep^{-(1+\alpha)}\,\|\hat{\cl
R_\ep}(\cdot,\cdot)\|_{L^1([0,1]^2)}= 0,\quad
\lim_{\ep\to 0} \ep^{-(1+\alpha)}\,\|\tilde{\cl
R_\ep}(\cdot;T)\|_{L^1([0,1])}= 0.
$$
Then, Assumption \ref{ass-ref2} holds, with
$$
\bar \varphi(t,s)=\int_0^t du \int_0^s dv \, \hat
g(u,v)\quad\mbox{and}\quad \bar \psi(t,T)=\int_0^t du \, \tilde
g(u;T).
$$

\end{proposition}

The proof is straightforward and postponed to Appendix B.

\begin{example}\label{mIBM-Y}\rm \bf[$m$-fold integrated
Brownian motion] \rm
Let us come back to Example \ref{IGP-Y-deg}, with $X$ as
the $m$-fold integrated Brownian motion:
$$
X_t=\int_0^tdu\Big(\int_0^{u}du_{m-1}\cdots\int_0^{u_2} du_1
W_{u_1}\Big),
$$
where $W$ denotes a standard Brownian motion.
Recall that the covariance function  is here
$$
k(t,s)=\int_0^t\int_0^s\kappa(u,v)\,du\,dv,\quad\mbox{with }
\kappa(t,s)=\frac 1{((m-1)!)^2}\int_0^{t\wedge s}
(t-\xi)^{m-1}(s-\xi)^{m-1}\,d\xi.
$$
Things are slightly different  according to $m=1$ or $m\geq 2$.
Let us consider  only  $m\geq 2$, the case $m=1$
being contained in next Example \ref{IFBM-Y}.
Straightforward computations allow to show that
\begin{equation}\label{app-1}
\begin{array}{c}
\kappa(T_n+\ep u, T_n+\ep v)=
\kappa(T_n,T_n)+\ep\frac 1{2}T_n^{2m-2}(u+v)+\smallskip\\
\displaystyle
+\ep^2 \frac {(m-1)}{2m-3}T_n^{2m-3}\Big[\frac{(m-2)}2(u+v)^2
+uv\Big]+\cl O(\ep^3)
\end{array}
\end{equation}
and
\begin{equation}\label{app-2}
\begin{array}{c}
\displaystyle
\int_0^T \kappa(T_n+\ep u,v) \, dv=
\int_0^T \kappa(T_n,v) \, dv+ \ep u
\frac{m-1}m \int_0^T (T_n-x)^{m-2}(T-x)^m\, dx+\smallskip\\
\displaystyle
+ \ep^2 u^2
\frac{(m-2)(m-1)}m\int_0^T (T_n-x)^{m-3}(T-x)^m\, dx+\cl
O(\ep^3),
\end{array}
\end{equation}
in which $\cl O(\ep^3)$ denotes a function going to $0$ as
$\ep\to 0$ in the right $L^1$ space at speed $\ep^3$.
Therefore, thanks to (\ref{app-1}), (\ref{app-2}) and
Proposition   \ref{prop-IGPref12}, Assumption
\ref{ass-ref2} does hold with
$$
\begin{array}{l}
\displaystyle
\bar\varphi(t,s)=
\frac{(m-1)}{4(2m-3)}T_n^{2m-3}[(m-2)(t^2s+ts^2)
+t^2s^2]\smallskip\\
\displaystyle
\bar\psi(t,T)= \Big(\frac 12
\frac{(m-2)(m-1)}m\int_0^T (T_n-x)^{m-3}(T-x)^m\, dx\Big) t^2.
\end{array}
$$
By using
Theorem \ref{th-bridgeref2}, we can assert that
the bridge process $Y$ satisfies a (non degenerate)
large deviation principle with inverse speed $\ep^4$ and
asymptotic covariance as in (\ref{barkyref2}), with $\alpha=1$.
\end{example}

\begin{example}\rm\label{IFBM-Y} \bf[Integrated
fractional Brownian motion] \rm Let $X_t=\int_0^tZ_udu$, with $Z$  a
fractional Brownian motion with Hurst index $H$. This is a quite
interesting example because, according to $H\leq 1/2$ or
$H>1/2$, one has both the cases studied in Proposition
\ref{prop-IGPref12}. In fact, straightforward computations
allow to state that
\begin{equation}\label{app-3}
\kappa_H(T_n+\ep u,T_n+\ep v)=
\kappa_H(T_n,T_n)+
HT_n^{2H-1} (u+v)\ep -\frac 12|u-v|^{2H}\ep^{2H} +\cl O(\ep^2)
\end{equation}
and
\begin{equation}\label{app-4}
\begin{array}{c}
\displaystyle\int_0^T \kappa_H(T_n+\ep u, v)dv=\\
=\left\{
\begin{array}{ll}
\displaystyle\int_0^T \kappa_H(T_n,
v)dv+\ep\Big[(HT_n^{2H-1} T -\frac 12 T_n^{2H} +\frac
12(T_n-T)^{2H})u\Big]\!+\!\cl O(\ep),& \!\!\! H<\frac 12\smallskip\\
T^2/2, &\!\!\! H=\frac 12\smallskip\\
\displaystyle\int_0^T \kappa_H(T_n, v)dv+\ep\Big[(HT_n^{2H-1} T -\frac 12
T_n^{2H} +\frac 12(T_n-T)^{2H})u\Big]\!+\!\cl O(\ep^2),&\!\!\! H>\frac
12\end{array}\right.
\end{array}
\end{equation}

Therefore, the asymptotic behavior can be resumed as follows.
\begin{itemize}
\item[(a)] If $H< 1/2$, part 1. in Proposition
\ref{prop-IGPref12} holds, with
$\alpha=2H$ and
$$
\hat g(u,v;T_n)=-\frac 12
|u-v|^{2H}\qquad \tilde g(u;T_n,T)\equiv 0.
$$

\item[(b)] If $H=1/2$, again part 1. in Proposition
\ref{prop-IGPref12} holds, with
$\alpha=1$ and
$$
\hat g(u,v;T_n)=\frac
12((u+v)-|u-v|)=u\wedge v \qquad \tilde g(u;T_n,T)\equiv 0.
$$

\item[(c)] If $H>1/2$, part 2. in Proposition
\ref{prop-IGPref12} holds, with
$\alpha=2H-1$ and
$$
\hat g(u,v;T_n)=-\frac 12 |u-v|^{2H} \qquad \tilde g(u;T_n,T)=0
$$

\end{itemize}

In conclusion, by suitably applying  Theorem \ref{th-bridgeref1} and
\ref{th-bridgeref2}, the family of bridges
$(Y^n_{T_n+\ep \cdot})_\ep$ satisfies a large
deviation principle on $C([0,1])$, with inverse speed
$\ep^{2+2H}$ and asymptotic  covariance function given by
$$
\bar k_Y(t,s)=\bar \varphi_H(t,s)+ t s\, \bar \varphi_H(1,1)- t\,
\bar \varphi_H(1,s)-s \,\bar \varphi_H(t,1),
$$
where
$$
\bar\varphi_H(t,s)\equiv \bar\varphi_n(t,s)=\left\{
\begin{array}{ll}\displaystyle
\frac {(|t-s|^{2H+2}-t^{2H+2}-s^{2H+2})}{2 (2H+1)(2H+2)}  & H\neq 1/2\smallskip\\
\displaystyle\frac{(t\wedge s)^3}3+\frac{(t\wedge s)^2}2|t-s| & H=
1/2.\end{array}\right.$$

Let us add some further remarks. In the case $H=1/2$, it is immediate
to check that
$\bar\varphi_H(t,s)= \int_0^t\int_0^s\kappa_H(u,v)\,dudv$. In other words,
$\bar\varphi_{1/2}$ turns out to be the covariance function
of the process $X$. Then, by taking into account Remark \ref{false-bridge},
the large deviations associated
to the bridge of the $n$-fold integrated Brownian motion
behave as ``the false bridge'', even if with a faster speed
(in fact, in this case the inverse speed is $\ep^3$, while the inverse
speed of the non-conditioned $n$-fold process is given by $\ep^2$).

\end{example}

\section{The asymptotic behavior of the crossing probability}
\label{exit}

In this section, the previous results
are applied in order to state the large deviation
asymptotic behavior of the hitting probability, the underlying
process of interest being the bridge of an
$n$-fold conditional Gaussian process. The already collected results can
be resumed in the following

\begin{hypothesis}\label{hyp}\rm
\begin{enumerate}
\item The family of the $n$-fold conditional processes
$\{(X^n_{T_n+\ep
t})_{t\in[0,1]}\}_\ep$ satisfies a large deviation principle with
inverse speed $\gamma^2_\ep$  and rate function
$$
J_n(h)=\left\{
\begin{array}{cl}
\displaystyle \frac 12\, \|h-x_n\|^2_{\bar{\cl H}_n}
& \mbox{ if } h_0=x_n \mbox{ and } h-\bar m\in \bar{\cl H}_n\\
+\infty & \mbox{ otherwise}
\end{array}
\right.
$$
$\bar{\cl H}_n$ being the reproducing kernel Hilbert space
associated to a suitable covariance function $\bar k_n$.

\item
The family of the bridges of the $n$-fold conditional processes
$\{(Y^n_{T_n+\ep
t})_{t\in[0,1]}\}_\ep$ satisfies a large deviation principle with
inverse speed $\eta^2_\ep$  and rate function
$$
J_Y(h)=\left\{
\begin{array}{cl}
\displaystyle \frac 12\, \|h-\bar m\|^2_{\bar{\cl H}_Y}
& \mbox{ if } h_0=x_n, h_1=y\mbox{ and } h-\bar m\in \bar{\cl H}_Y\\
+\infty & \mbox{ otherwise}
\end{array}
\right.
$$
$\bar{\cl H}_Y$ being the reproducing kernel Hilbert space
associated to a suitable covariance function $\bar k_Y$ and $\bar m_t
=x_n+\bar\beta_t(y-x_n)\equiv x_n+\bar k_n(t,1)(y-x_n)/\bar k_n(1,1)$.
\end{enumerate}
\end{hypothesis}

Throughout this section, we assume that Hypothesis \ref{hyp}
always holds.

\smallskip

Now, let us first focus on the upper
barrier case, the same arguments will apply for lower barriers.

Let $U:\R\rightarrow \R$ be a continuous function standing for an
upper barrier, and consider the probability
that $Y^n_{T_n+\ep \cdot}$ reaches the barrier $U$ up to the final time 1,
that is
$$
\P(\tau^U_\ep\leq 1),\quad
\mbox{with } \tau^U_\ep=\inf\{t>0\,:\,Y^n_{T_n+\ep t}\geq U_{T_n+\ep t}\}
$$
The above probability is negligible if $Y^n_{T_n}=x_n<U_{T_n}$ and
$Y^n_{T_n+\ep t}=y<U_{T_n+\ep t}$ for any $\ep$ close to 0, that
is $y\leq U_{T_n}$. As we will see, the case $y= U_{T_n}$ will
give a non relevant estimate, so that we can assume both $x_n$ and
$y$ are less than $U_{T_n}$.
So, if $x_n, y<U_{T_n}$, one has
$$
\eta^2_\ep\log \P(\tau^U_\ep\leq 1)\cong -I_Y^U,
$$
as $\ep \cong 0$, with $I_Y^U>0$.
Let us now see what  $I_Y^U$ is. Set $Z^n_{T_n+\ep t}=Y^n_{T_n+\ep
t}-U_{T_n+\ep t}$.
Since $\lim_{\ep\to 0} U_{T_n+\ep t}=U_{T_n}$ uniformly for
$t\in[0,1]$, by contraction it immediately follows that
$\{(Z^n_{T_n+\ep
t})_{t\in[0,1]}\}_\ep$ satisfies a large deviation principle as well,
with the same inverse speed and rate function
$$
J_Z(h)=J_Y(h+U_{T_n}).
$$
Then, one has
$$
\displaylines{
\lim_{\ep\to 0}\eta^2_\ep\log\P(\tau^U_\ep\leq 1)=-\inf_{\gamma\in\Gamma_U}J_{Y}(\gamma+U_{T_n})=-I_Y^U ,
}
$$
being $\Gamma_U=\{\gamma\,:\,
\sup_{t\in[0,1]} \gamma_t\geq 0\}$.

If a (continuous) lower barrier $L_t$ were considered, then the same arguments
would apply, giving
$$
\lim_{\ep\to 0}\eta^2_\ep\log\P(\tau^L_\ep\leq 1)
=-\inf_{\gamma\in\Gamma_L}J_{Y}(\gamma+L_{T_n})=-I_Y^L ,
$$
where $\tau^L_\ep=\inf\{t>0\,:\,Y^n_{T_n+\ep t}\leq L_{T_n+\ep t}\}$
and $\Gamma_L=\{\gamma\,:\, \inf_{t\in[0,1]} \gamma_t\leq 0\}$,
and this is interesting when $x_n, y>L_{T_n}$.
Finally, in the double barrier case, with $L_t\leq U_t$ for any $t$, then
the hitting probability behaves as follows:
$$
\lim_{\ep\to 0}\eta^2_\ep\log\P(\tau^{L,U}_\ep\leq 1)
=-I_Y^{L,U},
$$
where $\tau^{L,U}_\ep=\tau^L_\ep\wedge \tau^U_\ep$ is the first time
at which $Y^n_{T_n+\ep \cdot}$ reaches at least one barrier
and $I_Y^{L,U}$ is a suitable quantity, which is strictly positive
if $x_n,y\in (L_{T_n},U_{T_n})$.

The quantities $I_Y^U$, $I_Y^L$ and $I_Y^{L,U}$ are computed in next

\begin{proposition}\label{I-UL}
Suppose that $L$ and $U$ are continuous functions, with $L_t\leq U_t$
for any $t\in[0,1]$. Then,
$$
\begin{array}{lll}
I_Y^{U}&=\displaystyle\inf_{t\in[0,1]} \frac {\Big((U_{T_n}-x_n)(1-\bar \beta_t)+
\bar\beta_t(U_{T_n}-y)\Big)^2}{2\,\bar k_Y(t,t)}
& \mbox{if }x_n, y<U_{T_n}\smallskip\\
I_Y^{L}&=\displaystyle\inf_{t\in[0,1]} \frac {\Big((x_n-L_{T_n})(1-\bar \beta_t)+
\bar\beta_t(y-L_{T_n})\Big)^2}{2\,\bar k_Y(t,t)}
& \mbox{if }x_n, y>L_{T_n}\smallskip\\
I_Y^{L,U}&=\min\Big(I_Y^{L},I_Y^{U}\Big)
& \mbox{if }x_n, y\in(L_{T_n}, U_{T_n})
\end{array}
$$

\end{proposition}

\proof Consider the first equality (sigle upper
barrier case). We have to show that
$$
\inf_{\gamma\in\hat\Gamma_U}\frac 12\,\|\gamma+U_{T_n}-\bar m\|^2_{\bar{\cl H}_Y}
=\inf_{t\in[0,1]} \frac {\Big((U_{T_n}-x_n)(1-\bar \beta_t)+
\bar\beta_t(U_{T_n}-y)\Big)^2}{2\,\bar k_Y(t,t)},
$$
being $\hat\Gamma_U=\{\gamma: \gamma+U_{T_n}-\bar m\in \bar{\cl H}_{Y}, \,
\sup_{t\in[0,1]} \gamma_t\geq 0\}$.
Setting
$\hat\Gamma_{t,U}=\{\gamma: \gamma+U_{T_n}-\bar m\in \bar{\cl H}_{Y}, \,
\gamma_t= 0\}$, one has that $\hat\Gamma=\cup_{0<t< 1}\hat\Gamma_t$, so
that we simply need that
$$
\inf_{\gamma\in\hat\Gamma_{t,U}}
\frac 12\,\|\gamma+U_{T_n}-\bar m\|^2_{\bar{\cl H}_Y}
=\frac {\Big((U_{T_n}-x_n)(1-\bar \beta_t)+
\bar\beta_t(U_{T_n}-y)\Big)^2}{2\,\bar k_Y(t,t)},
$$
As already seen (see Section \ref{tutorial}),
a set of paths which is dense in ${\bar{\cl H}}_Y$ is
the one formed by those which are the barycenters of the random variable
belonging to the dual space of $C([0,1])$, that is,
$$
\gamma_u+U_{T_n}-\bar m_u=\int_0^1 \bar k_Y(u,v)\lambda(dv),
$$
as $\lambda$ varies in $\cl M[0,1]$. Since for such kind of paths
$$
\|\gamma+U_{T_n}-\bar m\|^2_{\bar{\cl H}_Y}
=\int_0^1\int_0^1 \bar k_Y(u,v)
\lambda(du)\lambda(dv),
$$
it is enough to minimize the r.h.s. above with respect to $\lambda$,
with the additional constraint $\gamma_t=0$, which becomes
$$
\bar m_t-U_{T_n}+
\int_0^1 \bar k_Y(t,v)
\lambda(dv)=0.
$$
This is a constrained extremum problem: using Lagrange multipliers,
$\lambda$ must satisfy
$$
\int_0^1 \bar k_Y(u,v)\lambda(dv)-\alpha \bar k_Y(t,u)=0,
\quad \mbox{for any } u\in [0,1],
$$
for some $\alpha\in \R$. Taking care of the constraint one finds
$$
\alpha =\frac {U_{T_n}-\bar m_t}{\bar k_Y(t,t)},
\quad \quad \lambda(dv)=\frac
{U_{T_n}-\bar m_t}{\bar k_Y(t,t)}\,\delta_{\{t\}}(dv).
$$
$\delta_{\{t\}}$ standing for the Dirac mass in $t$. Therefore,
$$
\inf_{w\in \hat\Gamma_{t,U}} \frac 12 \int_0^1 \int_0^1
\bar k_Y(u,v)\lambda(u)\lambda(dv)= \frac
{(U_{T_n}-\bar m_t)^2}{2\bar k_Y(t,t)}
$$
and the statement immediately follows by recalling that
$\bar m_t=x_n+\bar\beta_t(y-x_n)$.

Concerning the second equality, it follows by developing analogous
arguments. As for the final one, it is standard in large
deviation theory (see e.g. the discussion in the proof of Theorem
2.2 in Baldi and Caramellino \cite{bib:bp}).
\cvd

Let us stress that the barriers $U$ and/or $L$ can be also
piecewise continuous, in which case
the previous machinery runs again if the jump times coincide
with some of the conditional times $T_1, \ldots, T_n$.

\smallskip

Before to develop some examples, let us recall that
$\bar\beta_t=\bar k_n(t,1)/\bar k_n(1,1)$, where
$\bar k_n$ is defined in (\ref{kbar-n}) and represents
 the asymptotic covariance function associated to
the $n$-fold conditional process.
When our first set of large deviation results for the bridge holds
(as in Section \ref{immediate-LD}), one has
$$
\bar k_Y(t,s)=\bar k_n(t,s)-\frac{\bar k_n(t,1)\bar k_n(s,1)}
{\bar k_n(1,1)},
$$
$\bar k_Y$ being more complicated if it turns out following  Section \ref{refined-LD}.
Therefore, the minimization problem as required to compute $I_Y^U$ and
$I_Y^L$ has not a closed form in general, so that for practical
purposes one might be forced to use some numerical method (e.g. the Newton
method).

\begin{example}\label{IFBM-Y-exit}\rm
\bf[Integrated Brownian motion] \rm
Following Example \ref{IFBM-Y} (with $H=1/2$), let us consider
$X_t=\int_0^t B_s\,ds$, with $B$ as a standard Brownian motion.
Such a process has interesting applications in metrology, where it
is used as a model for the atomic clock error and the exit from a
fixed boundary means that the clock error exceeds an allowed limit
so that it has to be re-synchronized. Here,
we are in the second set of our large deviation estimates:
the bridge of the $n$-fold conditional process satisfies
a large deviation principle at inverse speed $\ep^3$ and
with rate function associated to the asymptotic covariance function
$$
\begin{array}{l}
\displaystyle\bar k_Y(t,s)
=\bar \varphi(t,s)+ t s\, \bar \varphi(1,1)- t\,
\bar \varphi(1,s)-s \,\bar \varphi(t,1), \smallskip\\
\mbox{with}\quad \displaystyle\bar \varphi(t,s)=
\frac{(t\wedge s)^3}3+\frac{(t\wedge s)^2}2\,|t-s|.
\end{array}
$$
Since $\bar k_n(t,s)=a^2_n\,ts$, one has $\bar\beta_t=t$,
so that $I^U_Y=g(U_{T_n})$ and $I^L_Y=g(L_{T_n})$, with
$$
g(a)=\inf_{t\in[0,1]} \frac {\Big(|a-x_n|(1-t)+
t|a-y|\Big)^2}{2\,t^2(1-t)^2/3}.
$$
The solution is simple to find:
$$
g(a)=\frac 32\,\Big(|a-x_n|^{1/2}+|a-y|^{1/2}\Big)^4.
$$
\end{example}

\begin{example}\label{FBM-Y-exit} \rm
\bf[Fractional Brownian motion] \rm
Following Section \ref{FBM} and Example \ref{FBM-Y}, let us consider
a fractional Brownian motion $X$ with Hurst index $H$,
in which one has
$$
\bar k_n(t,s)=k_H(t,s)=\frac{t^{2H}+s^{2H}-|t-s|^{2H}}2.
$$
So, in order to compute $I^U_Y$ and $I^L_Y$, giving the asymptotic
behavior of the hitting probability of the bridge $Y^n$,
by Proposition \ref{I-UL} one should be able to compute
$$
g_H(a)=\inf_{t\in[0,1]} \frac {\Big((a-x_n)(1-k_H(t,1))+
k_H(t,1)(a-y)\Big)^2}{2\,\Big(k_H(t,t)-k^2_H(t,1)\Big)},
$$
either with $a>x_n,y$ or $a<x_n,y$. In fact, one has $I^U_Y=g_H(U_{T_n})$
and $I^L_Y=g_H(L_{T_n})$. As far as we know, the exact solution can
be computed only when $H=1/2$, that is when a standard
Brownian motion is taken into account, in which case one has
$$
g_{1/2}(a)=2(a-x_n)(a-y),
$$
which agrees with well known formulas
(see e.g. Baldi and Caramellino \cite{bib:bc}).

\end{example}

In relation to the above Example \ref{FBM-Y-exit},
we have performed some numerical experiments concerning the
fractional Brownian motion. In particular, we have estimated via
Monte Carlo the probability of crossing the upper barrier $U=1$ up to time $1$
in two different ways: by crude simulations, in which the exit is
reached if a simulated position is greater than $U=1$, and
by means of the corrected procedure as recalled in the Introduction,
for which at each step the crossing is decided by using the
large deviation approximation for the exit probability of the pinned
process. In all the experiments, the exit probability is numerically computed
through ${10}^5$ simulations.
The results are given in terms of the method (corrected/crude) and
of the step size
($\ep=0.01, 0.002, 0.001$), for varying
values of the Hurst index $H$, which is set equal to $0.3$, $0.5$ and $0.7$.
Whenever $H=0.5$, everything is known (exit probability$=0.31732$), including
the fact that the crude approach works very poorly, so it has ben considered
to asses the procedure and for comparison purposes.
The choices $H=0.3$ and $H=0.7$ have been taken
to compare the results when $H<1/2$ (short memory, more
irregular paths) and $H>1/2$ (long memory, less irregular paths).
The results, given in Table \ref{Tab:exit-probability-fbm-030507},
show how much is the sensitivity w.r.t. the method (corrected/crude)
when $H$ decreases, that is when the irregularity of the path tends
to be higher. This is not surprising because the
inverse speed of the large deviations for the bridge is in fact $\ep^{2H}$,
so that the correction works more when $H$ decreases.

{\small
\begin{table}[tb]
\begin{center}
\begin{tabular}{llccc}
${\phantom{\Big(}}$ {\sc Method} & {\sc Step} & $H=0.3$ & $H=0.5$ & $H=0.7$ \\
\hline
${\phantom{\Big(}}$ corrected &  $0.01$
&
$\begin{array}{c}
0.60876\\
\mbox{{\footnotesize $( 0.60573 , 0.61178 )$}}
\end{array}$
&
$\begin{array}{c}
0.31820\\
\mbox{{\footnotesize $( 0.31531 , 0.32109 )$}}
\end{array}$
&
$\begin{array}{c}
0.20564\\
\mbox{{\footnotesize $( 0.20313 , 0.20814 )$}}
\end{array}$
 \\
\hline
${\phantom{\Big(}}$ corrected & $0.002$ &
$\begin{array}{c}
0.61841\\
\mbox{{\footnotesize $ ( 0.61540 , 0.62142 )$}}
\end{array}$
&
$\begin{array}{c}
0.31980 \\
\mbox{{\footnotesize $ ( 0.31691 , 0.32269 )$}}
\end{array}$
&$\begin{array}{c}
0.20274 \\
\mbox{{\footnotesize $ ( 0.20025 , 0.20523 )$}}
\end{array}$
\\
\hline
${\phantom{\Big(}}$ crude & $0.01$
&$\begin{array}{c}
0.47909 \\
\mbox{{\footnotesize $ ( 0.47599 , 0.48219 )$}}
\end{array}$
&$\begin{array}{c}
0.28918 \\
\mbox{{\footnotesize $ ( 0.28637 , 0.29199 )$}}
\end{array}$
&$\begin{array}{c}
0.19884\\
\mbox{{\footnotesize $ ( 0.19637 , 0.20131 )$}}
\end{array}$
 \\
\hline
${\phantom{\Big(}}$ crude & $0.002$
&$\begin{array}{c}
 0.54114  \\
\mbox{{\footnotesize $ ( 0.53805 , 0.544230 )$}}
\end{array}$
&$\begin{array}{c}
 0.30496 \\
\mbox{{\footnotesize $ ( 0.30211 , 0.30781 )$}}
\end{array}$
&$\begin{array}{c}
0.20222\\
\mbox{{\footnotesize $ ( 0.19973 , 0.20471 )$}}
\end{array}$
\\
\hline
${\phantom{\Big(}}$ crude  & $0.001$
&$\begin{array}{c}
0.56082\\
\mbox{{\footnotesize $ ( 0.55774 , 0.56390 )$}}
\end{array}$
& $\begin{array}{c}
0.30878\\
\mbox{{\footnotesize $ ( 0.30592 , 0.31164 )$}}
\end{array}$
&$\begin{array}{c}
0.20251\\
\mbox{{\footnotesize $ ( 0.20002 , 0.20500 )$}}
\end{array}$
\\
\hline
\end{tabular}
\caption{Fractional Brownian motion: Monte Carlo estimated
probability of crossing the upper barrier $U=1$ up to time $1$,
for varying values of the Hurst index $H$. In brackets,
the associated 95\% confidence interval. }
\label{Tab:exit-probability-fbm-030507}
\end{center}
\end{table}
}


\hyphenation{pre-fe-ren-ze}

\section*{Appendix A: proof of (\ref{J-bridge-1})}
This Appendix is devoted to the proof of representation (\ref{J-bridge-1}).
Let $J_Y$ be the rate function given by Theorem \ref{th-bridge},
i.e.
$$
J_Y(h)=\left\{
\begin{array}{cl}
\displaystyle \frac 12\, \|h-\bar m\|^2_{\bar{\cl H}_Y}
& \mbox{ if } h_0=x_n, h_1=y, h-\bar m\in \bar{\cl H}_Y\\
+\infty & \mbox{ otherwise}
\end{array}
\right.
$$
where $\bar m_t =x_n+\bar\beta_t(y-x_n)$ and
$\bar{\cl H}_Y$ is the reproducing kernel Hilbert space
associated to the covariance function
$$
\bar k_Y(t,s)=\bar k_n(t,s)-\bar\beta_s\bar k_n(t,1) =\bar
k_n(t,s)-\frac{\bar k_n(t,1)\bar k_n(s,1)}{\bar k_n(1,1)},
$$
$\bar k_n$ being defined in
(\ref{kbar-n}). Then, $J_Y$ can be written as
$$
J_Y(h)=\left\{
\begin{array}{cl}
\displaystyle \frac 12\,\Big( \|h-x_n\|^2_{\bar{\cl
H}_n}-\frac{(y-x_n)^2}{\bar k_n(1,1)}\Big)
& \mbox{ if } h_0=x_n, h_1=y, h-x_n\in \bar{\cl H}_n\\
+\infty & \mbox{ otherwise}
\end{array}
\right.
$$
$\bar{\cl H}_n$ being the reproducing kernel Hilbert space
associated to the covariance function $\bar k_n$.

Let us observe that this can be done in two ways: by
large deviation arguments (in particular, by using contraction
type properties allowing to transfer large deviation
principles)  or by handling reproducing kernel
Hilbert spaces. Here, we follow the second way.

\smallskip

 First, let us prove that the sets where the two
functionals are finite are the same, that is
$\cl K_1=\cl K_2$, being
$$
\displaylines{
\cl K_1=\{h\,:\,h_0=x_n, h_1=y, h-\bar m\in \bar{\cl H}_Y\}\cr
\cl K_2=\{h\,:\,h_0=x_n, h_1=y, h-\bar m\in \bar{\cl H}_n\}.
}
$$
If we set
$$
\begin{array}{rl}
\cl D_1&:=\displaystyle \{h\in \cl K_1\,:\, h_t-\bar m_t=\int_0^1
\bar k_Y(t,s)\alpha(ds), \mbox{ for some } \alpha\in \cl M[0,1]\}
\smallskip\\
\cl D_2&:=\displaystyle \{h\in \cl K_2\,:\, h_t-\bar m_t=\int_0^1
\bar k_n(t,s)\gamma(ds), \mbox{ for some } \gamma\in \cl M[0,1]\},
\end{array}
$$
then the statement is a consequence of the fact that
\begin{equation}\label{appoggio}
\cl D_1=\cl D_2\quad\mbox{and}\quad \|h-\bar m\|_{\bar{\cl H}_Y}
=\|h-\bar m\|_{\bar{\cl H}_n},\mbox{ for any }h\in \cl D_1=\cl
D_2.
\end{equation}
Indeed, since $\overline{\cl D_1}^{\|\cdot \|_{\bar{\cl H}_Y}}=\cl
K_1$ and $\overline{\cl D_2}^{\|\cdot \|_{\bar{\cl H}_n}}=\cl
K_2$, it immediately will follow that $\cl K_1=\cl  K_2$. So, let
us show that (\ref{appoggio}) does hold.

If one takes $h\in \cl D_1$, then
$$
\begin{array}{rl}
h_t-\bar m_t=&\displaystyle\int_0^1 \bar k_Y(t,s)\alpha(ds)
=\int_0^1 \Big(\bar k_n(t,s)
-\frac{\bar k_n(1,t)\bar k_n(1,s)}{\bar k_n(1,1)}\Big)\alpha(ds)\smallskip\\
=&\displaystyle\int_0^1 \bar k_n(t,s)\Big(\alpha(ds)
-\frac{\int_0^1\bar k_n(1,u)\alpha(du)}{\bar
k_n(1,1)}\delta_{\{1\}}(ds)\Big)
\end{array}
$$
where $\delta_{\{1\}}$ denotes the Dirac mass, and then $h\in\cl
D_2$. Conversely, if $h\in\cl D_2$,
then $h_t-\bar m_t=\int_0^1\bar k_n(t,s)\gamma(ds)$, and
in particular it must be
$$
0=h_1-\bar m_1=\int_0^1 \bar k_n(1,s)\gamma(ds).
$$
Therefore,
$$
h_t-\bar m_t=\int_0^1 \bar k_n(t,s)\gamma(ds) =\int_0^1 \Big(\bar
k_n(t,s)
-\frac{\bar k_n(1,t)\bar k_n(1,s)}{\bar k_n(1,1)}\Big)\gamma(ds)\smallskip\\
=\int_0^1 \bar k_Y(t,s)\gamma(ds)
$$
and $h\in \cl D_1$. Finally,
$$
\|h-\bar m\|^2_{\bar{\cl H}_Y}=\int_0^1\int_0^1 \bar
k_Y(t,s)\alpha(ds)\alpha(dt) =\int_0^1\int_0^1 \bar
k_n(t,s)\gamma(ds)\gamma(dt) =\|h-\bar m\|^2_{\bar{\cl H}_n}
$$
where $\alpha$ and $\gamma$ denote the measures representing
$h-\bar m$ in $\cl D_1$ and $\cl D_2$ respectively, so that
(\ref{appoggio}) is completely proved.

Now, we need to prove that for any $h\in \cl D_2$ one has
$\|h-\bar m\|^2_{\bar{\cl H}_n} =\|h-x_n\|^2_{\bar{\cl
H}_n}-(y-x_n)^2/\bar k_n(1,1)$. This follows from the fact that
$\bar m-x_n$ belongs to the reproducing kernel Hilbert space
$\bar{\cl H}_n$, because
$$
\bar m_t-x_n=\frac{\bar k_n(t,1)}{\bar k_n(1,1)}\, (y-x_n)
=\int_0^1\bar k_n(t,s)\,\frac{y-x_n}{\bar
k_n(1,1)}\,\delta_{\{1\}}(ds).
$$
Moreover, it holds
$$
\|\bar m-x_n\|^2_{\bar{\cl H}_n} =\frac{(y-x_n)^2}{\bar
k_n(1,1)}.
$$
Take now $h\in \cl D_2$. In particular, for some
 measure $\gamma$ one has $h_t-\bar m_t=\int_0^1\bar k_n(t,s)\gamma(ds)$.
Then, the measure $\hat\gamma(ds)=\gamma(ds)+(y-x_n)\delta_{\{1\}}(ds)/\bar k_n(1,1)$
is such that $h_t-x_n=\int_0^1 \bar k_n(t,s)\,\hat\gamma(ds)$
and
$$
\displaylines{
\<h-x_n, \bar m-x_n\>_{\bar{\cl H}_n} =\int_0^1\int_0^1\bar
k_n(t,s)\Big(\gamma(ds)+\frac{y-x_n}{\bar
k_n(1,1)}\,\delta_{\{1\}}(ds)\Big)\,\frac{y-x_n}{\bar
k_n(1,1)}\,\delta_{\{1\}}(dt)\cr
=\frac{y-x_n}{\bar
k_n(1,1)}\int_0^1\int_0^1\bar
k_n(1,s)\gamma(ds)
=\frac{y-x_n}{\bar
k_n(1,1)}(h-x_n)_1=
\frac{(y-x_n)^2}{\bar
k_n(1,1)}\cr
 }
$$

Therefore,
$$
\|h-\bar m\|^2_{\bar{\cl H}_n} =\|h-x_n\|^2_{\bar{\cl
H}_n}+\|\bar m-x_n\|^2_{\bar{\cl H}_n} -2\<h-x_n, \bar m-x_n\>_{\bar{\cl H}_n} =\|h-x_n\|^2_{\bar{\cl
H}_n}-\frac{(y-x_n)^2}{\bar k_n(1,1)},
$$
and the statement finally holds. \cvd

\section*{Appendix B: proof of Proposition \ref{prop-IGPref12} }

\proof  1.  Since $X$ is an integrated Gaussian process, one has
$$\displaylines{\Cov(X_{T_n+\ep t}-X_{T_n}, X_{T_n+\ep
s}-X_{T_n})=\int_{T_n}^{T_n+\ep t} du\int_{T_n}^{T_n+\ep t} dv
\kappa(u,v)=\cr \ep^2 \int_0^t du \int_0^s dv\,
 \kappa(T_n+\ep u, T_n+\ep v) .}
$$
Therefore,  one has
$$\displaylines{ \Cov(X_{T_n+\ep t}-X_{T_n}, X_{T_n+\ep
s}-X_{T_n})=\ep^2 \int_0^t du \int_0^s  dv\Big[\kappa(T_n,T_n)
+\ep^\alpha \hat g(u,v) + \hat{\cl R}_\ep(u,v)\Big]=\cr
=\ep^2\Big(\kappa(T_n,T_n) \, ts + \ep^\alpha \int_0^t du\int_0^s
dv \,\hat g(u,v)+\int_0^t du\int_0^s dv\, \hat{\cl
R}_\ep(u,v)\Big),}$$ so that $(i)$ of  Assumption
\ref{ass-ref1} is satisfied with $\bar \varphi(t,s)=\int_0^t du
\int_0^s dv \, \hat g(u,v).$

Moreover, since
$$\displaylines{k(T_n+\ep t, T)-k(T_n
,T)=\int_{T_n}^{T_n+\ep t} du \int_0^T dv \,\kappa(u,v)= \ep
\int_0^t  du \int_0^T dv\, \kappa(T_n+\ep u, v), }
$$
one obtains
$$\displaylines{k(T_n+\ep t, T)-k(T_n,T)=\ep \Big(\int_0^t du \Big[\int_0^T dv\,\kappa(T_n,v) +
\ep^\alpha\tilde  g(u;T) + \tilde{\cl R}_\ep(u;T)\Big]\Big)=\cr
=\ep \Big( t\,\int_0^T dv\,\kappa(T_n,v)+\ep^\alpha\int_0^t du\,
\tilde g(u;T)+\int_0^t du\,\tilde{\cl R}_\ep(u;T)\Big).
 }$$
Then, also $(ii)$ of Assumption \ref{ass-ref1} is satisfied
with $\bar \psi(t,T)=\int_0^t du \, \tilde g(u;T).$

The proof of part 2. follows exactly the same lines as for part 1.
\cvd

\newpage

\fontsize{10pt}{12pt}
\section*{Appendix C}

\subsection*{C.1 Proof of (\ref{app-1})}
 One has,
$$
\displaylines{
\kappa(T_n+\ep u, T_n+\ep v)=\int_0^{T_n+\ep(u\wedge v)}
(T_n+\ep u-x)^{m-1}(T_n+\ep v-x)^{m-1} \,dx=\cr \sum_{k=0}^{m-1}
\sum_{j=0}^{m-1-k}
\coeffbin{m-1}{k}\coeffbin{m-1-k}{j}\ep^{2m-2-2k-j}(u+v)^j
(uv)^{m-1-k-j}\times\cr
\times \int_0^{T_n+\ep(u\wedge v)} (T_n-x)^{2k+j}\,dx=\cr
=\sum_{k=0}^{m-1} \sum_{j=0}^{m-1-k}
\coeffbin{m-1}{k}\coeffbin{m-1-k}{j}\ep^{2m-2-2k-j}(u+v)^j
(uv)^{m-1-k-j}\frac{T_n^{2k+j+1}}{2k+j+1} +\cr
-\ep^{2m-1}\sum_{k=0}^{m-1}
\sum_{j=0}^{m-1-k}(-1)^{2k+j+1}\coeffbin{m-1}{k}\coeffbin{m-1-k}{j}
(u+v)^j (uv)^{m-1-k-j}\frac{(u\wedge v)^{2k+j+1}}{2k+j+1}=\cr
=\frac {1}{2m-1}T_n^{2m-1}+ \frac 1{2}T_n^{2m-2}(u+v)\ep +\frac
{(m-1)}{2m-3}T_n^{2m-3}\Big[\frac{(m-2)}2(u+v)^2 +uv\Big]\ep^2
+\cl O(\ep^3)=\cr
=\kappa(T_n,T_n)+\frac 1{2}T_n^{2m-2}(u+v)\ep
+\frac {(m-1)}{2m-3}T_n^{2m-3}\Big[\frac{(m-2)}2(u+v)^2
+uv\Big]\ep^2 +\cl O(\ep^3).
}
$$
\subsection*{C.2 Proof of (\ref{app-2})}
$$\displaylines{ \int_0^T \kappa(T_n+\ep u,v) \, dv= \int_0^T dv\int_0^{T_n+\ep(u\wedge v)}
(T_n+\ep u-x)^{m-1}(v-x)^{m-1}\, dx=\cr
=\sum_{k=0}^{m-1} \coeffbin
{m-1}{k}\ep^{m-1-k}u^{m-1-k}\int_0^T \, dv\int_0^v (T_n-x)^k
(v-x)^{m-1} \, dx=\cr
=\int_0^T \kappa(T_n,v) \, dv+ \ep u (m-1)
\int_0^T dv\int_0^v (T_n-x)^{m-2} (v-x)^{m-1} \, dx +\cr
=\ep^2 u^2
\frac 12 (m-1)(m-2) \int_0^T dv \int_0^v (T_n-x)^{m-3} (T-x)^m \,
dx+ \cl O(\ep^3)=\cr
=\int_0^T \kappa(T_n,v) \, dv+ \ep u
\frac{m-1}m \int_0^T (T_n-x)^{m-2}(T-x)^m\, dx+ \cr
+\ep^2 u^2
\frac{(m-2)(m-1)}m\int_0^T (T_n-x)^{m-3}(T-x)^m\, dx+\cl
O(\ep^3).}
$$
\subsection*{C.3 Proof of (\ref{app-3})}
One has
$$
\displaylines{
\kappa(T_n+\ep u,T_n+\ep v)= \frac 12 [(T_n+\ep u)^{2H}+(T_n+\ep
v)^{2H}- |u-v|^{2H}\ep^{2H}]=\cr
=T_n^{2H} + HT_n^{2H-1} (u+v)\ep
-\frac 12|u-v|^{2H}\ep^{2H} +\cl O(\ep^2)=\cr
=\kappa(T_n,T_n)+
HT_n^{2H-1} (u+v)\ep -\frac 12|u-v|^{2H}\ep^{2H} +\cl O(\ep^2).
}
$$
\subsection*{C.4 Proof of (\ref{app-4})}
For $v\leq T_n$ and $T\leq T_n$, one has
$$
\displaylines{
\int_0^T \kappa(T_n+\ep u, v)dv=\frac 12 \int_0^T
((T_n+\ep u)^{2H}+v^{2H}-(T_n- v+\ep u)^{2H})\, dv=\cr
=\int_0^T
\kappa(T_n, v)dv+ \frac 12 \int_0^T dv\,((T_n+\ep
u)^{2H}-T_n^{2H})-\frac 12 \int_0^T dv\,((T_n- v+\ep u)^{2H}-(T_n-
v)^{2H}).
}
$$
Straightforward calculations show that
$$\displaylines{
\int_0^T \kappa(T_n+\ep u, v)dv=\cr
=\left\{
\begin{array}{ll}\int_0^T \kappa(T_n,
v)dv+\ep\Big[(HT_n^{2H-1} T -\frac 12 T_n^{2H} +\frac
12(T_n-T)^{2H})u\Big]+\cl O(\ep)& H<\frac 12\\
\frac{T^2}2 & H=\frac 12\\
\int_0^T \kappa(T_n, v)dv+\ep\Big[(HT_n^{2H-1} T -\frac 12
T_n^{2H} +\frac 12(T_n-T)^{2H})u\Big]+\cl O(\ep^2)& H>\frac
12\end{array}\right.
}
$$


\addcontentsline{toc}{section}{References}

\begin{thebibliography}{99}

\bibitem{bib:azencott}
R. Azencott (1980). Grande D\'eviations et applications. In
\'Ecole d'\'et\'e de probabilit\'es de St. Flour VIII, L.N.M. Vol
774, Springer, Berlin/Heidelberg/New York.

\bibitem{bib:bci}
P. Baldi and L. Caramellino, M.G. Iovino (1999). Pricing general
barrier options: a numerical approach using sharp large
deviations. \emph{Mathematical Finance}, \bf 9,  \rm 293-321.

\bibitem{bib:bc}
P. Baldi and L. Caramellino (2002). Asymptotics of hitting
probabilities for general one-dimensional pinned diffusions.
\emph{Annals of Applied Probability}, \bf 12\rm, 1071-1095.

\bibitem{bib:bp}
P. Baldi and B. Pacchiarotti (2006). Explicit computation of second
order moments of Importance Sampling
estimators for fractional Brownian motion.
\emph{Bernoulli}, \bf 12\rm, 663--688.

\bibitem{bib:cl}
X. Chen and W.V. Li (2003). Quadratic functionals and small ball
probabilities for the $m$-fold integrated Brownian motion.
\emph{Annals of Probability}, \bf 31\rm, 1052-1077.

\bibitem{bib:cheridito1}
P. Cheridito (2001).
Mixed fractional Brownian motion. \emph{Bernoulli}, \bf 7\rm, 913-934.


\bibitem{bib:dz}
A. Dembo and O. Zeitouni (1992). \emph{Large deviations techniques
and applications}. Jones and Bartlett Publishers.

\bibitem{bib:ds}
J.D. Deuschel and D.W. Stroock (1989). \emph{Large deviations}.
Academic press, Boston.

\bibitem{bib:torino}
L. Galleani, L. Sacerdote, P. Tavella and C. Zucca (2003). A mathematical
model for the atomic clock error. \emph{Metrologia}, \bf 40\rm, 257-264.

\bibitem{bib:gsv}
D. Gasbarra, T. Sottinen and E. Valkeila (2004).
Gaussian bridges. Preprint available at \texttt{http://math.tkk.fi/reports/a481.pdf}.


\bibitem{bib:gjw}
P. Groeneboom, G. Jongbloed and J.A. Wellner (2001).
A canonical process for estimation of convex functions:
the ``invelope'' of integrated Brownian motion $+t\sp 4$.
\emph{Annals of Statistics}, \bf 29\rm, 1620-1652.


\bibitem{bib:mk}
G. Molchan and A. Khokhlov (2004). Small values
of the maximum or the integral of fractional Brownian motion.
\emph{Journal of Statistical Physics}, \bf 114\rm, 923-946.

\end{thebibliography}
\end{document}